\input amstex

\documentstyle{amsppt}
\magnification=\magstep1
\def\mn{\medskip\noindent}
\def\subheading#1{\medbreak\line{\hfil\bf#1\hfil}\smallskip}
\topmatter
\title Closed geodesics on orbifolds \endtitle
\author K. Guruprasad
and A. Haefliger \endauthor
\address{Universit\'e de Gen\`eve, Section de Math\'ematiques, 2-4 Rue du
Li\`evre, 1211 Gen\`eve
24}\endaddress

\thanks{Most of the results of this paper were obtained in 1990 while
 the first named author  was visiting
Geneva. On January 5. 1998, K. Guruprasad
 died in a tragic accident in Bangalore.}
\endthanks

\abstract{In this paper, we try to generalize to the case of compact
Riemannian orbifolds $Q$ some
classical results about the existence of closed geodesics of positive
length on compact
Riemannian manifolds $M$. We shall also consider the problem of the
existence of infinitely
many geometrically distinct closed geodesics.

In the classical case the solution of those problems involve the
consideration of the homotopy
groups of $M$ and the homology properties of the free loop space on $M$
(Morse theory). Those
notions have their analogue in the case of orbifolds (see  \cite{9}). The
main part of this paper
will be to recall those notions and to show how
the classical techniques can be adapted to the case of orbifolds.} \endabstract

\endtopmatter

\document

\heading 1. Summary of the results \endheading

 A Riemannian orbifold structure $Q$ on a Hausdorff
topological space $|Q|$ is given by an atlas of uniformizing charts $q_i:
X_i \to V_i$. The
$V_i$, $ i \in I$, are open sets whose union is $|Q|$; each $X_i$ is a
Riemannian manifold with a
finite subgroup $\Gamma_i$ of its group of isometries, the $q_i$ are continuous
$\Gamma_i$-invariant maps  inducing a homeomorphism from $\Gamma_i
\backslash X_i$ to $V_i$.
The change of charts are Riemannian isometries.

	It is convenient to assume that the $X_i$ are disjoint and to
consider on the union $X$ of the
$X_i$ the pseudogroup $\Cal P$ generated by the change of charts (its
restriction to $X_i$ is
generated by the elements of $\Gamma_i$), or equivalently the topological
groupoid
$\Cal G$ of germs of
change of charts. It is an \'etale groupoid with space of units $X$, and
the projection $q$
which is the union of the $q_i$ induces a homeomorphism from the space of
orbits $\Cal G
\backslash X$ to $|Q|$.

The orbifold $Q$ is said to be {\it developable} ({\it good} in the sense
of Thurston \cite{17}) if it is
the quotient of a Riemannian manifold by a discrete subgroup $\Gamma$ of
its group of
isometries.

A continuous free loop on $Q$ is an equivalence (cohomology) class of
continuous $\Cal
G$-cocycles (see 2.2.1) on the circle
$S^1$. To put a topology on the set $|\Lambda^c Q|$ of continuous  free
loops on $Q$, it is
convenient to represent such a free loop
 as an equivalence class (see 2.3) of closed
$\Cal G$-paths
$c=(g_0,c_1,g_1,\dots,c_k,g_k)$ over a subddivision
$0=t_0<t_1<\dots<t_k=1$ of the interval $[0,1]$. Here the $c_i$ are
continuous maps from
$[t_{i-1},t_i]$ to $X$, the  $g_i \in \Cal G$ are germs of changes of
charts with source
$c_{i+1}(t_i)$ and target $c_i(t_i)$, for $0<i<k$, the source of $g_0$
(resp. $g_k$) being
$c_1(0)$ (resp. the target of $g_0$) and the target of $g_k$ being
$c_k(1)$. The length of the
free loop
$[c]$ represented by $c$ is the sum of the length of the $c_i$, and if the
$c_i$ are
differentiable, its energy is  equal to the sum of the energies of the
$c_i$. It is a closed geodesic
on $Q$ if each
$c_i$ is a constant speed geodesic and the differential of $g_i$, for
$1<i<k$,  maps the velocity
vector $\dot c_{i+1}(t_i)$ to $\dot c_i(t_i)$ and the differential of
$g_kg_0$ maps $\dot
c_1(0)$ to $\dot c_k(1)$. Free loops of length $0$ are always geodesics.
They can be
represented by pairs of the form $(c_1,g_1)$, where $c_1:
[0,1] \to X$ is the constant map to a point $x$ and $g_1$ an element of
$\Cal G$ with source and
target $x$ (well defined up to conjugation).

A continuous $\Cal G$-loop based at $x \in X$ is also an equivalence class
of $\Cal
G$-paths $c$ as above,
where  $x$ is the  target of $g_0$ (which is egal to the source of $g_k$),
under a more
restrictive equivalence class (where the point
$x$ is preserved, see 2.3.2 for the precise definition). Let $\Omega^c_X =
\Omega^c_X(\Cal G)$ be the
set of continuous
$\Cal G$-loops based at the various points of $X$. The groupoid $\Cal G$
acts naturally on $\Omega^c_X$
and there is a natural projection
$q_\Lambda:
\Omega^c_X \to |\Lambda^c Q|$ inducing a bijection  $\Cal G \backslash
\Omega^c_X
\to |\Lambda^c
Q|$.

To each orbifold $Q$ is associated a classifying space $BQ$ (see \cite{9}).
Roughly
speaking, it is a space
$BQ$ with a projection
$\pi: BQ \to |Q|$ such that, for each $y \in |Q|$, the fiber $\pi^{-1}(y)$
is a space whose
universal covering in contractible and whose fundamental group is
isomorphic to the
isotropy subgroup of a point $x \in q_i^{-1}(y)$ in a uniformizing chart
$q_i: X_i \to V_i
\subseteq |Q|$. The elements of  $|\Lambda^c Q|$ are in bijection with
equivalence classes of
continuous maps $\Bbb S^1 \to BQ$,
two maps being  equivalent if they are homotopic through an homotopy
projecting under $\pi$ to a
constant homotopy. The fundamental group of $BQ$ is isomorphic to the
orbifold fondamental group
of $Q$.

For the proof of the following proposition, see 3.1 and 3.2.

\proclaim{ Proposition}  $\Omega^c_X$  is naturally a Banach manifold. The
action of $\Cal G$ is differentiable and the quotient space $|\Lambda^c Q|$
has a natural
structure of  Banach orbifold noted $\Lambda^c Q$.

	If $BQ$ and $B\Lambda^c Q$ denote the classifying spaces of the
orbifolds $Q$ and $\Lambda^c Q$,
then
$B\Lambda^c Q$ has the same weak homotopy type as the topological space
$\Lambda BQ$ of continuous
free loops on $BQ$. \endproclaim

The fact that the free loop space of an orbifold has an orbifold structure
was also observed independently
by Weimin Chen \cite{3}.

Let  $|\Lambda Q|$ (resp. $\Omega_X$) be the sets of $\Cal G$-loops (resp.
based $\Cal G$-loops) of class $H^1$,
namely those  loops represented by closed
$\Cal G$-paths
$c= (c_1,g_1,\dots,c_k,g_k)$ such that each path $c_i$ is absolutely continuous
and the velocity
 $|\dot c_i|$ is square integrable. The energy $E(c)$ of $c$
is equal to the sum of the energies of the $c_i$ and depends only on its
equivalence class. The next proposition is proved in 3.3.2 and 3.3.5.

\proclaim{ Proposition} $\Omega_X$ is naturally a Riemannian manifold
and $|\Lambda Q|$ has
a Riemannian orbifold structure noted $\Lambda Q$.

	The natural inclusion $\Lambda Q \to \Lambda^c Q$ induces a
homotopy equivalence on the
corresponding classifying spaces. \endproclaim

One proves like in the classical case that the energy function $E$ is
differentiable on $\Lambda Q$ and
that its critical points correspond to the closed geodesics. If the
orbifold $Q$ is compact, the
Palais-Smale condition holds. The usual techniques
(Lusternik-Schnierelmann, Morse) as explained for
instance in Klingenberg \cite{11} or \cite{12}, can then be applied to study
the existence of closed
geodesics on compact Riemannian orbifolds. We get in particular the
following results,
proved in 5.1.

\proclaim{ Theorem} Let $Q$ be a compact Riemannian orbifold. There is a
closed geodesic on
$Q$ of positive length in each of the following cases:

1) $Q$ is not developable,

2) The (orbifold) fundamental group of $Q$ has an element of infinite order
or is finite.

 \endproclaim

We don't know if this covers all the cases. It would be so if a finitely
presented torsion group
would be finite (see the remark 5.1.2).

We say that two closed geodesics on $Q$ are geometrically distinct if their
projections to
$|Q|$ have distinct images. The following result (see 5.2.3) is an
extension   to
orbifolds of a classical result of Gromoll-Meyer \cite{5} and
Vigu\'e-Sullivan \cite{18}.

\proclaim{ Theorem} Let $Q$ be a compact simply connected Riemannian orbifold
such that the rational cohomology of $|Q|$ is not generated by a single
element.
Then there exist on
$Q$ an infinity of geometrically distinct closed geodesics of positive
length.  \endproclaim

The paper is organized as follows.
 The lengthy section 2 recalls elementary and classical definitions
concerning orbifolds,
''maps'' from spaces to orbifolds, in particular closed curves (or free
loops) on orbifolds,
classifying spaces of orbifolds, etc. In section 3 we explain the orbifold
structure on the
space $|\Lambda^c Q|$ of continuous free loops and on the space $|\Lambda
Q|$ of free loops of class $H^1$
on $Q$ and the basic properties of the energy function. In section 4 we
study the orbifold
tubular neighbourhoods of the $\Bbb S^1$-orbits of closed geodesics in a
geometric equivalence
class. In the last section, we sketch how the
classical theory of closed geodesics on Riemannian manifolds can be adapted
to the case of
orbifolds. In sections 3, 4 and 5 we assume familiarity with the
  notions and the basic papers concerning  the theory of closed geodesics
 in classical Riemannian
geometry.

\medskip

A first version of this paper was submitted to Topology. In his report, the
referee pointed out
to us the rich literature (we were not aware of) concerning the theory of
invariant
closed geodesics, in particular the work of Grove-Tanaka \cite{6}. Thanks
to his remarks we
were able in this revised version to extend to compact simply connected
Riemannian orbifold the
result of Gromoll-Meyer
\cite{5}, using \cite{6}. The interpretation in our framework of the
results of Grove-Tanaka
and Grove-Halperin \cite{7} is mentionned in 3.3.6  and 5.2.5.

\heading 2. Basic definitions  \endheading

In this section we recall basic definitions and known results (see \cite{9}
and \cite{2}).

\subheading{2.1 Orbifolds}

\mn{\bf 2.1.1.  Definition of an orbifold structure.}
Let $|Q|$ be a Hausdorff topological space. A
differentiable orbifold structure
$Q$ on
$|Q|$ is given by the following data:

i) an open cover $\{V_i\}_{i \in I}$ of $|Q|$ indexed by a set $I$,

ii) for each $i \in I$, a finite subgroup $\Gamma_i$ of the group of
diffeomorphisms of a connected
differentiable manifold $X_i$ and a continuous  map $q_i:X_i \to V_i$,
called a uniformizing chart,
inducing a homeomorphism from
$\Gamma_i\backslash X_i$ onto $V_i$,

iii) for each $x_i \in X_i$ and $x_j \in X_j$ such that $q_i(x_i) =
q_j(x_j)$, there is a
diffeomorphim $h$ from an open connected neighbourhood $W$ of $x_i$ to a
neighbourhood of
$x_j$ such that $q_j \circ h = {q_i}_{|W}$. Such a map $h$ is called a {\it
change of chart} ; it
is well defined up to composition with an element of $\Gamma_j$; if $i=j$,
then $h$ is the
restriction of an element of $\Gamma_i$.

The family $(X_i,q_i)$ is called an {\it atlas of uniformizing charts}
defining the orbifold
structure $Q$. By definition two such atlases define the same orbifold
structure on $|Q|$ if,
put together, they satisfy the compatibility condition iii). It is easy to
check that the above definition
of orbifolds is equivalent to the definition of $V$-varieties introduced by
Satake (as explained in
\cite{10}) and to the definition of Thurston \cite{17}.

A {\it Riemannian orbifold} is an orbifold $Q$ defined by an atlas of
uniformizing charts such
that the $X_i$ are Riemannian manifolds and the change of charts are
Riemannian isometries. Note
that on any paracompact differentiable orbifold one can introduce a
Riemannian metric.

\mn{\bf 2.1.2. The pseudogroup of change of charts. Developability.}
Let
$X$ be the disjoint union of the
$X_i$. We identify each $X_i$ to a connected component of $X$ and we note
$q: X \to |Q|$ the union
of the maps $q_i$. Any diffeomorphim $h$ from an open subset $U$ of $X$ to
an open subset of
$X$ such that $q \circ h = q_{|U}$ will be called a change of charts. The
collection of
change of charts form a pseudogroup $\Cal P$ of local diffeomorphisms of
$X$, called the
pseudogroup of change of charts of the orbifold (with respect to the
uniformizing atlas
 $(X_i,q_i)$). Two points $x,y \in X$ are said to be in the same orbit of
$\Cal P$ if there
is an element $h \in \Cal P$ such that $h(x)=y$. This defines an
equivalence relation on
$X$ whose classes are called the orbits of $\Cal P$. The quotient of $X$ by
this
equivalence relation, with the quotient topology, will be denoted $\Cal P
\backslash X$.
The map $q:X \to |Q|$ induces a homeomorphism from $\Cal P \backslash X$ to
$|Q|$.

The pseudogroups of change of charts of two atlases defining the
same orbifold structure on $|Q|$ are equivalent in the following sense. Two
pseudogroups
$\Cal  P_0$ and $\Cal P_1$ of local diffeomorphisms of differentiable
manifolds $X_0$ and
$X_1$ respectively are equivalent if there is a pseudogroup $\Cal  P$ of local
diffeomorphims of the disjoint union $X$ of $X_0$ and $X_1$ whose
restriction to $X_j$ is
equal to   $\Cal  P_j$ and such that the inclusion of $X_j$ into $X$ induces a
homeomorphism $\Cal  P_j \backslash X_j \to \Cal  P \backslash X$, $j =0,1$.

More generally, consider a pseudogroup of local diffeomorphisms $\Cal  P$ of a
differentiable manifold $X$ such that each point $x$ of $X$ has an open
neighbourhood $U$
such that the restriction of
$\Cal  P$ to $U$ is generated by a finite group $\Gamma_U$ of
diffeomorphisms of $U$.
Assume moreover that the space of orbits $\Cal  P \backslash X$ is
Hausdorff. Then $\Cal  P
\backslash X$ has a natural orbifold structure still noted $\Cal  P
\backslash X$.

	For instance if $\Gamma$ is a discrete subgroup of the group of
diffeomorphisms of a manifold
$X$ whose action on $X$ is proper, then $\Gamma \backslash X$ has a natural
orbifold
structure. An orbifold structure arising in this way is called {\it
developable}. In this case, the pseudogroup of change of charts is
generated by the
elements of $\Gamma$.

\mn{\bf 2.1.3. The tear drop (a non-developable orbifold) \cite{17}.} The
topological space $|Q|$ is the
$2$-sphere $\Bbb S^2$ with north pole $N$ and south pole $S$. The north
pole is a conical point of order
$n$, all the  other points are regular. The orbifold structure can be
defined by an atlas of two
uniformizing charts $q_i: X_i \to V_i$, $i=1,2$. Here $X_i$ is the open
disc of radius $3\pi/4$ centered at
$0$ in
$\Bbb R^2$ with polar coodinates $(r,\theta)$, the group $\Gamma_1$ is
generated by a rotation of order $n$
and
$\Gamma_2$ is trivial. The map
$q_1$ (resp. $q_2$) maps the point $(r,\theta)$ to the point of $\Bbb S^2$
with geodesic coordinates
$(r,n\theta)$ (resp. $(r,\theta)$) centered at  the north pole (resp. the
south pole). The map
$q_2^{-1}q_1$ is the $n$-fold covering $(r,\theta) \mapsto (\pi
-r,n\theta)$ of the annulus
$A:=\{(r,\theta),
\pi/4 < r< 3\pi/4\}$. On the disjoint union $X = (\{1\} \times X_1) \cup
(\{2\} \times X_2)$ of $X_1$ and
$X_2$, the pseudogroup of changes of charts is generated by $\Gamma_1$ and
by the
diffeomorphisms from open sets of $\{1\}
\times A$ to open sets of $\{2\} \times A$ which are restrictions of the
map $(1,a) \mapsto
(2,q_2^{-1}q_1(a))$. Any Riemannian metric  on $Q$ is given by a Riemannian
metric on $X_1$ invariant by
$\Gamma_1$ and a Riemannian metric on
$X_2$ such that $q_2^{-1}q_1$ is a local isometry.

\mn{\bf 2.1.4. The \'etale groupoid of germs of change of charts.}  Recall
that a groupoid $(\Cal
G,X)$ is a small category $\Cal G$ with set of objects $X$, all elements of
$\Cal G$ being
invertible. The set of objects $X$ is often identified to the set of units
of $\Cal G$ by
the map associating to an object $x \in X$  the unit $1_x \in \Cal G$. Each
element  $g
\in \Cal G$ is considered as an arrow with source $\alpha(g) \in X$
(identified to its right
unit)  and target
$\omega(g)
\in X$ (identified to its left unit). The inverse of $g$ is denoted
$g^{-1}$. For $x \in \Cal G$, the
elements $g$ of $\Cal G$ with $\alpha(g) = \omega(g)$ form a group noted
$\Cal G_x$, called the isotropy
subgroup of $x$.

A topological groupoid is a groupoid
$(\Cal G,X)$ such that $\Cal G$ and $X$ are topological spaces, all the
structure maps
(projections $\alpha,\ \omega: \Cal G \to X$, composition, passage to
inverse)  are continuous,
and such that the map
$x
\mapsto 1_x$ from
$X$ to
$\Cal G$ is a homeomorphism onto its image. An {\it \'etale groupoid} is a
topological groupoid such
that the projections $\alpha$ and $\omega$ from $\Cal G$ to $X$ are
\'etale, i.e. are local
homeomorphisms.

To the pseudogroup $\Cal P$ of change of charts of an  atlas of uniformizing
charts defining an orbifold, we can associate the
\'etale groupoid $(\Cal G,X)$ of all germs of change of charts, with the
usual topology of
germs, $X$ being the disjoint union of the sources of the charts. From
$\Cal G$ we can
reconstruct $\Cal P$, because its elements can be obtained as the
diffeomorphisms from
open sets $U$  of $X$ to open sets of $X$ which are the composition with
$\omega$ of
sections $U \to \Cal G$ of $\alpha$ above $U$. We shall also use the
notation $Q = \Cal
G \backslash X$ to denote an orbifold whose pseudogroup $\Cal P$ of change
of charts is equivalent to
the pseudogroup corresponding to $\Cal G$.

For instance, consider the case where $Q$ is developable, quotient of a
smooth manifold $X$ by a discrete group
$\Gamma$ of diffeomorphisms of $X$, acting properly on $X$. Then $\Cal G$
is the groupoid $\Gamma \times
X$, where $\Gamma$ is endoved with the discrete topology, the projection
$\alpha$ (resp. $\omega$) mapping
$(\gamma,x)$ to $x$ (resp. $\gamma.x$). The composition
$(\gamma,x)(\gamma',x')$ is defined when
$x=\gamma'.x'$ and is equal to $(\gamma\gamma',x')$. We shall use the
notation $\Cal G = \Gamma \ltimes X$ and the orbifold $Q = \Gamma
\backslash X$ is also
noted $\Cal G \backslash X$.

The following lemma will be useful later on.

\proclaim{2.1.5. Lemma} Assume that  $Q = \Cal G\backslash X$ is a
Riemannian orbifold, and let $g \in \Cal
G$ with
$x=
\alpha(g)$ and $y =\omega(g)$. Let $B(x,\epsilon)$ and $B(y,\epsilon)$ be
two convex geodesic balls centered
at $x$ and
$y$ with radius $\epsilon$ (such balls always exist for small enough
$\epsilon$). Then there is a unique element
$h$ of the pseudogroup of changes of charts $\Cal P$ which is an isometry
from $B(x,\epsilon)$ to
$B(y,\epsilon)$ and whose germ at $x$ is equal to $g$. \endproclaim

\demo{Proof} By hypothesis the exponential map is defined on the ball of
radius $\epsilon$ centered at the
origin of $T_xX$ and is a diffeomorphism of that ball to $B(x,\epsilon)$;
similarly for $y$. Therefore let
$h: B(x,\epsilon) \to B(y,\epsilon)$ be the diffeomorphism mapping
isometrically  geodesic rays issuing
from $x$ to those issuing from $y$ and whose differential at $x$ is the
differential $Dg$ of a representative of $g$. For $0 <r
<\epsilon$, assume that  the restriction $h_r$ of $h$ to the open ball
$B(x,r)$ belongs to $\Cal P$; this
is the case if $r$ is small enough. It will be sufficient to prove that the
restriction of $h$ to a
neighbourhood of the closure of $B(x,r)$ belongs to $\Cal P$. For a point
$z \in \partial B(x,r)$, the points
$z$ and
$h(z)$ are in the same orbit under
$\Cal P$, because
$|Q|$ is Hausdorff. There are small  ball neighbourhoods $U$ of $z$ and $V$
of $h(z)$  and  an element $f:U
\to V$ of $\Cal P$  such that the restriction of $\Cal P$ to $U$ is
generated by a group $\Gamma_U$ of
diffeomorphisms of $U$ and such that the germ
 of any element of $\Cal P$ with source in $U$ and target in $V$ is the
germ of the composition of $f$
with an element of $\Gamma_U$. Therefore, the restriction of $h_r$ to
$B(x,r) \cap U$ is the restriction
of an element of $\Cal P$ defined on $U$. As such an element is a
Riemannian isometry, it must
coincide with
$h_{|U}$ on the geodesic rays issuing from $x$, hence on $U$. $\square$
\enddemo

\subheading{2.2 Morphisms from spaces to topological groupoids}

\mn{\bf 2.2.1. Definition using cocycles \cite{8}.}
Let $\Cal  G$ be a topological  groupoid with space of units $X$, source
and target
projections $\alpha, \omega: \Cal  G \to X$ respectively. Let $\Cal  U =
(U_i)_{i \in I}$ be
an open cover of a topological space $K$. A $1$-cocycle over $\Cal  U$ with
value in $\Cal  G$ is a collection of continuous maps $f_{ij}: U_i \cap U_j
\to \Cal  G$ such that, for
each $x \in U_i \cap U_j \cap U_k$,  we have:
$$f_{ik}(x) = f_{ij}(x) f_{jk}(x).$$
This implies in particular that $f_{ii}(x)$ is a unit of $\Cal  G$ and that
$f_i:=f_{ii}$ can be
considered as a continuous map from $U_i$ to $X$. Also $f_{ij} = f^{-1}_{ji}$.

Two cocycles on two open covers of $K$ with value in $\Cal  G$ are
equivalent if there
is a cocycle with value in $\Cal  G$ on the disjoint union of  those two
covers extending
the given ones on each of them. An equivalence class of cocycles is called
a (continuous) morphism
from $K$ to $\Cal  G$ (or when $Q$ is an orbifold $\Cal  G \backslash
X$ a ''continuous map'' from $K$ to $Q$). The set of equivalence classes of
$1$-cocycles on
$K$ with value in $\Cal  G$ is noted $H^1(K,\Cal  G)$.
 One should observe that {\it if $\Cal G$ and $\Cal G'$ are the groupoids of
germs of the changes of
charts of two atlases defining the same orbifold structure on a space
$|Q|$, then there is a natural
bijection bvetween the sets $H^1(K,\Cal G)$ and $H^1(K,\Cal G')$}.

Any continuous map $f$ from a topological space $K'$ to $K$ induces a map
$f^*: H^1(K,\Cal
G)
\to H^1(K', \Cal G)$. Two morphisms from $K$  to $\Cal G$ are homotopic if
there
is a morphism from $K \times [0,1]$ to $\Cal G$ such that the morphisms
from $K$
to $\Cal G$ induced by the natural inclusions $k \mapsto (k,i)$, $i=0,1$,
from $K$ to $K
\times [0,1]$ are the given morphisms.

Any morphism from $K$ to $\Cal G$ projects, via $q: X \to \Cal G\backslash
X =|Q|$, to a
continuous map from $K$ to $|Q|$; note that two distinct morphisms may have
the same projection (see
the example at the end of 2.3.5).

\mn{\bf 2.2.2.  Principal $\Cal G$-bundles.}
Another description of morphisms from $K$ to $\Cal G$ can be given in terms
of isomorphism
classes of principal
 $\Cal  G$-bundles over $K$ (see \cite{9}).

Let $E$ be a topological space with a continuous map $\alpha_E: E \to X$.
Let $E
\times_X \Cal G$
 be the subspace of $E \times \Cal G$ consisting of pairs $(e,g)$ such that
$\alpha_E(e) = \omega(g)$.
A continuous (right) action of
$\Cal  G$ on $E$ with respect to $\alpha_E$ is a continuous map $(e,g) \mapsto
e.g$ from the space
$E \times_X \Cal G$ to $E$ such that
$\alpha(g) = \alpha_E(e.g)$, $(e.g).g' = e.(gg')$ and $e.1_x = e$. Left actions
are defined similarly.

 A  principal $\Cal G$-bundle over $K$ is a topological space $E$
together with a surjective continuous map $p_E:E \to K$, called the bundle
projection,  and a continuous
action
$(e,g)
\mapsto e.g$ of $\Cal  G$ on $E$
 with respect to a continuous map
$\alpha_E:E \to X$, called the action map, such that $p(e.g) = p(e)$.
Moreover we assume that the action
is simply transitive on the fibers of $p$ in the following sense. Each
point of $K$ has an open
neighbourhood
$U$ with a continuous section $s: U \to E$ with respect to $p_E$ such that
the map $U \times_X
\Cal G
\to p^{-1}(U)$ mapping pairs  $(u,g) \in U \times \Cal G$ with $\omega(g) =
\alpha_Es(u)$ to
$s(u).g$ is a homeomorphism. It follows that if $\Cal  U = (U_i)_{i \in I}$
is an open
cover of $K$ and if $s_i:U_i \to E$ is a local continuous section of $p$
above $U_i$ for
each
$i
\in I$,  then
there are unique continuous maps $f_{ij}:U_i \cap U_j \to \Cal  G$ such
that $s_i(u) =
s_j(u)f_{ji}(u)$ for each $u \in U_i \cap U_j$. Thus $f = (f_{ij})$ is a
$1$-cocycle over
$\Cal  U$ with value in $\Cal  G$.

Conversely, if $f = (f_{ij})$ is a $1$-cocycle over an open cover
$\Cal  U= (U_i)_{i \in I}$ of $K$ with value in $\Cal  G$, then we can
construct a
principal
 $\Cal  G$-bundle $E$ over $K$ by identifying in the disjoint union of the $U_i
\times_X \Cal G =
\{(u,g) \in U_i \times  \Cal G: \omega(g)= f_{ii}(u)\}$ the point $(u,g)
\in U_i
\times_X \Cal G$, $u \in U_i \cap U_j$,  with the point $(u,f_{ji}(u)g) \in
U_j \times_X
\Cal G$. The projections $p_E: E \to K$ and $\alpha_E: E \to X$ map the
equivalence class of $(u,g)
\in U_i \times_X \Cal G$ to $u$ and $\alpha(g)$ resp. and the action of
$g'$ on the class
of $(u,g)$ is the class of $(u,gg')$. A principal $\Cal G$-bundle obtained
in this way by
using an equivalent cocycle is isomorphic to the preceding one, i.e. there
is a
homeomorphism between them projecting to the
identity of $K$ and commuting with the action of $\Cal G$. This isomorphism
is determined
uniquely by a cocycle extending the two given cocycles.

Therefore we see that {\it there is a natural bijection between the set
$H^1(K,\Cal G)$ and the set of isomorphism classes of principal $\Cal
G$-bundles over
$K$}. This correspondence is functorial via pull back: if $E$ is a
principal $\Cal G$-bundle over
$K$ and if $f:K' \to K$ is a continuous map, then the pull back $f^*E$ of
$E$ by f (or the bundle
induced from $E$ by $f$) is the bundle $K' \times_K E$ whose elements are
the pairs $(k',e) \in K'
\times E$ such that $f(k') p_E(e)$. The projection (resp. the action map) sends
$(k',e)$ to $k'$ (resp. to $\alpha_E(e)$

$\Cal G$ itself can be considered as a principal $\Cal G$-bundle over $X$
with respect
to the projection $\omega: \Cal G \to X$, the action map $\alpha_\Cal G:
\Cal G \to X$
being the source
projection. Any principal $\Cal G$-bundle $E$ over $K$ is locally the pull back
of this bundle $\Cal G$ by a continuous map to $X$. The projection $p_E$ is an
\' etale map (i.e locally a homeomorphism) when $\Cal G$ is an \' etale
groupoid.

\medspace

 \mn{\bf 2.2.3.  Relative morphisms.}  Let $K$ be a topological
space,
$L \subseteq K$ be a subspace and $F$ be a principal $\Cal G$-bundle over
$L$. A morphism
from $K$ to $\Cal G$ relative to $F$ is represented by a pair $(E,\phi)$
where $E$ is a principal
$\Cal G$-bundles
$E$ over
$K$ and $\phi$ is  an isomorphism  from $F$ to the restriction $E_{|L}$ of $E$
above
$L$. Two such pairs $(E,\phi)$ and $(E',\phi')$ represent the same morphism
from $K$ to $\Cal
G$ relative to $F$ if there is an isomorphism $\Phi: E \to E'$ such that
${\phi'} = \Phi \circ
\phi$.

An important particular case is when $L$ is a base point $z \in K$ and $F$ is
the pull back of $\Cal G$ by the map sending $z$ to a given point $x \in X$,
i.e. $F$ is the set of pairs $(z,g)$, where $g$ is an element of $\Cal G$ with
target $x$.
  A morphism
from $K$ to $\Cal G$ relative to
$F$ (in that case we shall say a morphism from $K$ to $\Cal G$ mapping
$z\in K$ to $x \in X$) is just
represented  by a
$\Cal G$-bundle $E$ over $K$ with a base point $e \in E$ with $p_E(e) = z$ and
$\alpha_E(e)=x$. Indeed an isomorphism from $F$ to $p_E^{-1}(z)$ is determined
by the image $e$ of $(z,1_x) \in F$. Two such pairs
$(E,e)$ and $(E',e')$ represent the same morphism if there is an
isomorphism from
$E$ to $E'$ mapping $e$ to $e'$
and projecting to the identity of $K$. More generally, suppose that  $F$ is the
pull back of
$\Cal G$ by a continuous map $f_0:L \to X$. Then a morphism from $K$ to
$\Cal G$
relative to $F$ (we shall say relative to $f_0$) is
given by a principal
$\Cal G$-bundle
$E$ over
$K$ with a continuous  map $s: L \to E$ with $p_E \circ s$ the identity of
$L$ and $\alpha_E \circ s = f_0$.

Two morphisms represented by $(E_0,\phi_0)$ and $(E_1,\phi_1)$ from $K$ to
$\Cal G$
relative to $F$ are homotopic (relative to $F$) if there is a bundle $E$
over $K \times I$
and an isomorphism from $E_{|(K \times \partial I) \cup (L \times I)}$ to
the bundle obtained
by gluing  $F \times I$ to  $E_0 \times \{0\}$  and $E_1 \times \{1\}$
using the isomorphims
$\phi_0$ and $\phi_1$.

We leave to the reader a description of  isomorphism classes of relative
bundles in terms
of equivalence classes of relative cocycles.

Let $I^n = [0,1]^n$ be the $n$-cube, and let $\partial I^n$ be its
boundary. Fix a base
point $x$ in $X$. Let $f_0$ be the constant map from $\partial I^n $ to
$X$. We define
$\pi_n((\Cal G,X),x)$
as the set of homotopy classes of principal $\Cal G$-bundle over $I^n$
relative to $f_0$. One
proves as usual that this set has a natural group structure, called the
$n^{th}$- homotopy
group of $(\Cal G,X)$ based at $x$. In the case where $\Cal G \backslash X$
is a connected  orbifold $Q$,
this group is called the $n$-th homotopy group of $Q$, and for n=1 the
(orbifold) fundamental group of $Q$.

\subheading{2.3 Paths and loops in orbifolds}

In this section we describe in a more concrete
way the morphisms from the interval $I = [0,1]$ to a topological groupoid
$(\Cal G,X)$ relative
to a map from
$\partial I$ to $X$. Whe shall deal with
 an orbifold struture $Q$ on a topological space $|Q|$ defined by an atlas of
uniformizing charts,  $(\Cal G,X)$ will be the \'etale topological groupoid
of germs of
change of charts, or more generally the groupoid of germs of  elements of a
pseudogroup
defining $Q$ (cf. 2.1.2 and 2.1.4). As above, we have the map
$q:X
\to |Q|$ inducing a homeomorphism from the space of orbits $\Cal G
\backslash X$ to $|Q|$.
\thickspace

\mn { \bf 2.3.1. Continuous $\Cal G$-paths.} Let $x$ and $y$ be
two points of
$X$. A (continuous)
$\Cal G$-path from
$x$ to
$y$ over a subdivision
$0 =t_0 < t_1 < \dots < t_k=1$ of the interval $[0,1]$ is
a sequence
$c = (g_0,c_1,g_1,\dots, c_k,g_k)$ where

i) $c_i:[t_{i-1},t_i] \to X$ is a continuous map,

ii) $g_i$ is an element of $\Cal G$ such that $\alpha(g_i) =c_{i+1}(t_i)$ for
$i=0,1,\dots,k-1$, $\omega(g_i) = c_i(t_i)$ for $i=1,\dots,k$, and
$\omega(g_0) = x$, $ \alpha(g_k) = y$.

\thickspace

\mn { \bf 2.3.2. Equivalence classes of $\Cal G$-paths. The set
$\Omega^c_{x,y}$}.  Among
$\Cal G$-paths from
$x$ to
$y$ parametrized by
$[0,1]$ we define an equivalence relation generated by the following two
operations:

i) Given a $\Cal G$-path $c = (g_0,c_1,g_1,\dots, c_k,g_k)$ over the
subdivision $0=t_0 <
\dots < t_k=1$, we can add a subdivision point $t' \in (t_{i-1},t_i)$
together with the
unit element $g' = 1_{c_i(t')}$ to get a new sequence, replacing $c_i$ in
$c$ by
$c'_i,g',c''_i$, where $c'_i$ and $c''_i$ are the restriction of $c_i$ to the
intervals $[t_{i-1},t']$ and $[t',t_i]$.

ii) Replace the $\Cal G$-path $c$ by a new path $c' = (g'_0,c'_1,g'_1,\dots,
c'_k,g'_k)$ over the same subdivision as follows: for each $i=1,\dots,k$,
choose
continuous maps $h_i: [t_{i-1},t_i]\to  \Cal G$ such that $\alpha(h_i(t))
=c_i(t)$, and
define $c'_i: t \mapsto \omega(h_i(t))$, $g'_i =
h_i(t_i)g_i{h_{i+1}}(t_i)^{-1}$ for
$i=1,\dots,k-1$, $g'_0=g_0h_1(0)^{-1}$ and $g_k'= h_k(1)g_k$.

 The  equivalence class of a $\Cal
G$-path $c$ from $x$ to $y$ will be noted $[c]_{x,y}$, and the set
of such equivalence classes will be noted
$\Omega^c_{x,y}(\Cal G)$, or simply
$\Omega^c_{x,y}$ (here $c$ stands for continuous).  It corresponds
bijectively to the
set of isomorphims classes of principal $\Cal
G$-bundles $E$ over $I = [0,1]$ with two base points $e_0$ and $e_1$ over
$0$ and
$1$ such that $\alpha_E(e_0) = x$ and $\alpha_E(e_1) = y$ (see 2.2.3).  The
bundle
$E$ is obtained from $c$ as the quotient of the union of the bundles
$c^*_i(\Cal G)$ by the
equivalence relation identifying $(t_i, g_ig) \in c^*_i(\Cal G)$ to
$(t_i,g) \in
c^*_{i+1}(\Cal G)$ for $i=1,\dots ,k-1$. The base point $e_0$ is represented by
$(0,g_0) \in c_1^* (\Cal G)$ and $e_1$ by $(1,g_k^{-1})$.

If $c =(g_0,c_1,g_1,\dots,c_k,g_k)$ and $c'
=(g'_0,c'_1,g'_1,\dots,c'_k,g'_k)$ are two equivalent $\Cal G$-paths from
$x$ to $y$ over the same
subdivision, then the maps
$h_i$ in 2.3.2, ii) above are unique, because $\Cal G$ is Hausdorff and
\'etale. Therefore
there is a unique isomorphism from the relative principal $\Cal G$-bundle
associated to
$c$ to the one associated to $c'$.

Let   $c
=(g_0,c_1,g_1,\dots,c_k,g_k)$ be a $\Cal G$-path from $x$ to $y$ and $g$ be an
element of $\Cal G$ with $\alpha(g) = x$ and $\omega(g)= z$ (resp.
$\alpha(g) = z$
and $\omega(g)= y$). Then $g.c:=(gg_0,c_1,g_1,\dots,c_k,g_k)$ (resp.
$c.g:=(g_0,c_1,g_1,\dots,c_k,g_kg)$ ) is a $\Cal G$-path from $z$ to $y$ (resp.
from
$x$ to $z$) whose equivalence class depends only on the equivalence class
of $c$
and which is noted $g.[c]_{x,y} \in \Omega^c_{z,y}$ (resp $[c]_{x,y}.g \in
\Omega^c_{x,z}$).

\mn{\bf 2.3.3. Based $\Cal G$-loops. The sets $\Omega^c_x$ and $\Omega^c_X$}

A continuous $\Cal G$-path $c$ from $x$ to $x$ is called a closed $\Cal G$-path
(based at
$x$). Its equivalence class is called a continuous
$\Cal G$-loop based at $x$ and is noted $[c]_x$. The set of continuous
$\Cal G$-loops
based at $x$ is also noted
$\Omega^c_x(\Cal G)$ or simply $\Omega^c_x$.
We note $\Omega^c_X(\Cal G)$ or simply $\Omega^c_X$ the set  $\bigcup_{x \in
X}
\Omega^c_x$  of based continuous $\Cal G$-loops.

The set $\Omega^c_x$ is in bijection with the set of isomorphisms classes of
principal
$\Cal G$-bundle $E$  over the circle $\Bbb S^1= \Bbb R/\Bbb Z$ with a base
point
$e$ over
$0 \in \Bbb S^1$ such that
$\alpha_E(e) = x$ . A pointed principal $\Cal G$-bundle $(E,e)$
corresponding to
$c$ is constructed as follows. For
$i=1,\dots,k$,  let $E_i= [t_{i-1},t_i] \times_X \Cal G$ be the pull
back of the principal
$\Cal G$-bundle $\Cal G$ by the map $c_i:[t_{i-1},t_i] \to X$, i.e. the
space  of
pairs $(t,g)$ with $c_i(t)= \omega(g)$;  recall that the action map
$\alpha_i:E_i \to
X$ sends
$(t,g)$ to
$\alpha(g)$. The bundle
$E$ is the quotient of the
disjoint union of the
$E_i$ by the equivalence relation identifying $(t_i, g) \in
E_{i+1}$ to $(t_i, g_ig)
\in E_i$ for
$1 \le i \le k$ and $(0,g) \in E_1$ to $(1, g_kg_0g) \in E_k$. The base
point $e$
is the equivalence class of
$(0,g_0^{-1}) \in E_1$, or equivalently of $(1,g_k) \in E_k$. The only
automorphism of a pointed principal $\Cal G$-bundle $(E,e)$ over $\Bbb S^1$
(respecting the base point and projecting to the identity of $\Bbb S^1$) is the
identity.

The groupoid
$\Cal G$ acts naturally on the left on the set $\Omega^c_X$ of  based
continuous  $\Cal
G$-loops with respect to the
projection $ \Omega_X \to X$ associating to a  $\Cal G$-loop based at $x$ the
point $x$. Indeed let   $c =(g_0,c_1,g_1,\dots,c_k,g_k)$ be a closed $\Cal
G$-path
based at
$x$ and
$g$ be an element of $\Cal G$ with $\alpha(g) = x$ and $\omega(g)= y$. The
action of
$g$ on the equivalence class of $c$ is the $\Cal G$-loop $^g[c]_x$ based at $y$
represented by
$^gc:=  (gg_0,c_1,g_1,\dots,c_k,g_kg^{-1})$.

 If $[c]_x$ is represented
by the pointed
principal
$\Cal G$-bundle
$(E,e)$,  then $^g[c]_x$ is represented by the pointed bundle  $(E,
e.g^{-1})$.

For an integer $m$, the map $t \mapsto mt \mod 1$ from
$\Bbb S^1$ to $\Bbb S^1$
induces a map $\Omega^c_x \to \Omega^c_x$ mapping $[c]_x$ to $[c]^m_x$. On the
level of pointed bundles, it is just the pull back. Concretely when $m$ is a
positive integer, if
$c=(g_0,c_1,\dots, c_k,g_k)$ is a closed $\Cal G$-path over a subdivision
$0=t_0< \dots < t_k=1$, the image of $[c]_x$ by this map is equal to $[c^m]_x$,
where $c^m$ is the m-th
iterate of $c$. More precisely, $c^m= (g'_0,c'_1,\dots, g'_{mk})$ is the
closed $\Cal G$-path over
the subdivision $0=t'_0 < t'_1 < \dots < t'_{mk}=1$, where, for
$r=0,...,m-1$ and $i=1,\dots,k$, we
have
$t'_{rk+i} =\frac rm + \frac{t_i}m$, $c'_{rk+i}(t) = c_i(mt-r)$; we have
$g'_{(r+1)k}= g_kg_0$,
$g'_0 =g_0$ and $g'_{mk} = g_k$, and for $i \neq k$ we have
$g'_{rk+i}=g_i$.

\mn{\bf  2.3.4. The set $|\Lambda^c Q|$ of continuous free loops.} It is the
quotient $\Cal G\backslash \Omega^c_X$ of $\Omega_X$ by the action of $\Cal G$
described above. Its elements are called continuous free loops on $Q$.
They
correspond to the elements of
$H^1(\Bbb S^1,
\Cal G)$, the set of continuous morphisms from $\Bbb S^1$ to $\Cal G$. The
notation
is justified by the observation that
$H^1(\Bbb S^1,\Cal G) = H^1(\Bbb S^1,\Cal G')$ if $\Cal G$ and $\Cal G'$
are the
groupoids of germs of changes of charts of two atlases of uniformizing charts
defining the same orbifold structure on
$|Q|$.

An element of $|\Lambda^c Q|$  is
represented by a closed
$\Cal G$-path
$c = (g_0,c_1,g_1,
\dots, c_k,g_k)$ over a subdivision $0=t_0 < t_1 < \dots < t_k=1$ as in
2.3.2 based at some  point
$x
\in X$. This time the equivalence relation is generated   by i) and ii) in
2.3.2
and also by

iii) for any element $g \in \Cal G$ such that $\alpha(g)=x$, then $c =
(g_0,c_1,g_1,\dots,c_k,g_k)$ is equivalent to $^gc :=
(gg_0,c_1,g_1,\dots,c_k,g_kg^{-1})$.
The class of $c$ under this equivalence relation is noted $[c]$

Under the projection $q: X \to |Q|$, every free $\Cal G$-loop is mapped to
a free loop on $|Q|$.
Therefore if $\Lambda|Q|$ is the space of continuous free loops on the
topological
space $|Q|$ in the usual sense, we have a map
$$ |\Lambda^c Q| \to \Lambda |Q|.$$ This map is not injective in general
(see below).

 We note
$|\Lambda^0 Q|$ the subset of
$|\Lambda^c Q|$ formed by the free loops on $Q$ projecting to a constant
loop. An element of this
subset is represented by a closed $\Cal G$-path $c=(g_0,c_1,g_1)$, where
$g_0$ is a unit $1_x$, $c_1$
is the constant map from $[0,1]$ to $x$ and $g_1$ is an element of the
subgroup $\Cal G_x= \{ \{g
\in \Cal G: \alpha (g) = \omega (g) = x\}$. The equivalence class  $[c]$ of
$c$ correspond to the
conjugacy class of $g_1$ in $\Cal G_x$.

Alternatively $|\Lambda^c Q|$ can be described as the set of isomorphism
classes of principal $\Cal G$-bundle
$E$ over
$\Bbb S^1$. The group of homeomorphisms of $\Bbb S^1$ acts
on $|\Lambda^c Q|$ by change of parametrisation: if $h$ is a homeomorphism of
$\Bbb S^1$, its action on the
isomorphism class of a bundle $E$ is the isomorphism class of the pull back
of $E$ by $h$. In particular
the group $\Bbb S^1$ of rotations of $\Bbb S^1$ acts
on
$|\Lambda^c Q|$. The fixed point set of this action is precisely
$|\Lambda^0 Q|$.

\thickspace

\mn{\bf 2.3.5. The developable case.} Let $Q$ be the orbifold quotient of a
connected manifold X by the
action of a discrete subgroup
$\Gamma$ of its group of diffeomorphisms acting properly and let $\Cal G$
be the groupoid $\Gamma \ltimes
X$ (see the end of 2.1.4).  The set $\Omega^c_X$ of based continuous $\Cal
G$-loops  are in
bijection with the pairs $(c,\gamma)$, where $c:[0,1] \to X$ is a
continuous path and
 $\gamma$ is an element of $\Gamma$ mapping $c(0)$ to $c(1)$.
Indeed, consider a  $\Cal G$-loop at $x$
represented by the closed $\Cal G$-path $(g_o,c_1,\dots, c_k,g_k)$ over the
subdivision
$0=t_0 <t_1<\dots<t_k=1$, where $g_i = (\gamma_i,c_{i+1}(t_i))$ for
$i=0,\dots,k-1$ and
$g_k=(\gamma_k,c_k(1))$. Then we define
$\gamma =
\gamma_0
\dots
\gamma_k$ and $c:[0,1] \to X$ as the path defined for $t \in [t_{i-1},t_i]$
by $\gamma_0\dots \gamma_{i-1}.c_i(t)$.
 The continuous
free loops on
$Q$ are represented by classes of pairs $(c,\gamma)$ like above, such a pair
being equivalent to $(\delta
\circ c, \delta\gamma\delta^{-1})$, where $\delta \in \Gamma$. So
$|\Lambda^c Q|$ is the quotient of
$\Omega^c_X$ by this action of $\Gamma$.  Assuming
$X$ simply connected, the
set of homotopy classes  of elements of
$|\Lambda^c Q|$ is in bijection with the set of conjugacy classes in $\Gamma$.

We can equivalently describe the elements of $\Omega^c_X$ as the pairs
$(c,\gamma)$ where $c: \Bbb R \to X$ is
a continuous map such that
$c(t+1) =
\gamma.c(t)$. With this interpretation, we can describe the action of
$\Bbb S^1$ on $|\Lambda^c Q| =
\Gamma \backslash\Omega^c_X$ as follows. We have a natural action of
$\Bbb R$ on
$\Omega^c_X$ by translations: the action of  $\tau \in \Bbb R$
is given by
$\tau.(c,\gamma)=(c_\tau,\gamma)$, where
$c_\tau(t) = c(t + \tau)$.  This action commutes
 with the  action
$^\delta (c,\gamma):= (\delta\circ c,
\delta\gamma\delta^{-1})$ of
$\delta \in \Gamma$ on
$\Omega_X$. Therefore we get an action of $\Bbb R$ on the quotient
$\Gamma \backslash \Omega^c_X$. As
the translations by the integers act trivially, we get the action of $\Bbb
S^1 =
\Bbb R/\Bbb Z$ on
$|\Lambda^c Q|$.

As an example, let $X = \Bbb R^2$ and $\Gamma$ be the group generated by a
rotation $\rho$ fixing $0$
and of angle $2\pi/n$. Let $\Cal G = \Gamma\ltimes X$ be the groupoid
associated to
the action of $\Gamma$ on $X$. The
orbifold $Q = \Cal G \backslash X$ is a cone.  Consider the free
$\Cal G$-loop represented by the
pair
$(c,\rho^k)$, where $c$ is the constant path at $0$. If we deform this loop
slightly so that it
avoids the origin, its projection to the cone $|Q|$ will be a curve going
around the
vertex a number of times congruent to $k$ modulo $n$ ; in particular, when
$k=n$, it could also be a constant loop.

\subheading{2.4 Geodesics on Riemannian orbifolds}

\mn{\bf 2.4.1. Length and Energy.} We consider a Riemannian orbifold $Q =
\Cal G \backslash X$. The
length $L(c)$ of a $\Cal G$-path $c = (g_0,c_1,g_1,\dots, c_k,g_k)$ is the
sum of the length of
the paths $c_i$. It depends only on the equivalence class of $c$. If $c$ is
piecewise
differentiable (i.e. if each $c_i$ is piecewise differentiable), the length
of $c$ is given by
$$L(c) = \sum_1^k \int_{t_{i-1}}^{t_i} |\dot{c_i}(t)|\  dt,$$ and the
energy $E(c)$ of $c$
is
  $$E(c) = 1/2 \sum_1^k \int_{t_{i-1}}^{t_i} |\dot{c_i}^2(t)|\  dt.$$
>From Schwartz inequality we have
$$L(c)^2 \leq 2 E(c),$$
with equality if and only if the speed $|\dot{c_i}(t)|$ is constant for all
$i$ and $t$. The
length and the energy of $c$ depends only on its equivalence class.

The distance $d(z,z')$  between two points $z = q(x)$ and $z'= q(x')$ in
$|Q|$ is  defined as
the infimum of the length of $\Cal G$-paths joining $x$ to $x'$ (this is
independent of the
choice of $x$ and $x'$ in their $\Cal G$-orbit). This distance defines a
metric on $|Q|$.
The Riemannian orbifold $Q$ is
said to be complete, if this metric is complete. This is always the case if
$|Q|$ is compact.

\thickspace

 \mn{\bf 2.4.2. Geodesic $\Cal G$-paths and closed geodesics on orbifolds.}

 A geodesic $\Cal
G$-path from
$x$  to
$y$ in a Riemannian orbifold is a
$\Cal G$-path
$c = (g_0,c_1,g_1,\dots, c_k,g_k)$  from
$x$ to
$y$ such that  each $c_i$ is a geodesic segment with
constant speed and such that  the differential $Dg_{i+1}$ of a
representative of $g_{i+1}$ at
$c_{i+1}(t_i)$, maps the velocity vector
$\dot{c}_{i+1}(t_i)$ to the velocity vector $\dot{c}_{i}(t_i)$. Note that the
image of a $\Cal G$-geodesic path under the projection to $|Q|$ is an arc which
in general is not locally length minimizing for the metric on $|Q|$ defined
above (see the example in 2.4.5 ).

Note that if $c$ is a geodesic $\Cal G$-path
from $x$ to $y$, then the vector $Dg_0(\dot c_1(0))$ is an invariant of the
equivalence class
$[c]_{x,y}$ and is called the initial vector of the $\Cal G$-geodesic $c$.

If $c$ is a closed $\Cal
G$-path, it represents a closed $\Cal G$-geodesic $[c]_x$ based at $x$  if
moreover the differential of
$g_kg_0$ maps the velocity vector $\dot{c_1}(0)$ to the vector $\dot{c_k}(1)$.
Its equivalence class $[c]$ is called a closed geodesic on $Q$. A
free loop of length
$0$ is always a closed geodesic.

\mn{\bf 2.4.3. Geometric equivalence of closed geodesics.} Two closed geodesic
on $Q$ are geometrically equivalent if their image under the projection to
$|Q|$
are the same. Otherwise they are called geometrically distinct. The elements of
the
$\Bbb S^1$-orbit of a closed geodesic
$[c]$ are all geometrically equivalent to $[c]$, as well as the multiples
$[c^m]$
for $m \ne 0$ (see 2.3.3).

A closed geodesic is called primitive if its length is the minimum of the
lengths
of the closed geodesics in its geometric equivalence class.

\mn{\bf 2.4.4. The exponential morphism.}
If $Q$ is a complete
orbifold, given a vector $\xi \in T_xX$, there is always a geodesic $\Cal
G$-path $c=
(g_0,c_1,g_1,\dots,c_k,g_k)$ over a subdivision of $[0,1]$ issuing from
$x$ with initial vector $\xi$. Any geodesic $\Cal
G$-path issuing from $x$  with initial vector $\xi$ is equivalent to a
$\Cal G$-path obtained from $c$ by
replacing $g_k$ by any element $g_kg$ where $g \in \Cal G$ is any element such
that
$\omega(g) =
\alpha(g_k)$.

In fact we have the analogue of the {\it exponential map} for complete
orbifolds, namely for each point $x
\in X$ there is a morphism $exp_x: T_xX \to Q$ defined by a principal $\Cal
G$-bundle $E_x$ over $ T_xX$.
The elements of
$E_x$ above $\xi \in T_xX$ are the equivalence classes of geodesic  $\Cal
G$-paths with initial vector
$\xi$. The right action of $\Cal G$ on $E_x$ is as described above. There
is a canonical base point $e_x$
above
$0 \in T_x$ characterized by the following property. Let $s$ be the local
section
 of
$E_x$ defined on a small ball $U \subset T_xX$ with center $0$ mapping $0$ to
$e$; then its composition with
$\alpha_{E_x}$ is the usual exponential map $U \to X$.

	Given $g \in \Cal G$ with $\alpha(g) = x$ and $\omega(g) = y$ there
is a unique isomorphism $(E_x,e_x) \to
(E_y,e_y)$ of pointed principal $\Cal G$-bundle projecting to the
differential $Dg:T_xX \to T_yX$ of $g$.

\mn{\bf 2.4.5. The developable case.} If the orbifold $Q$ is the quotient
$\Gamma \backslash X$ of a
connected Riemannian manifold
$X$ by a discrete subgroup $\Gamma$ of its group of isometries, then  any
closed geodesic  on
$Q$ is  represented by a pair $(c,\gamma)$, where
$c:\Bbb R \to X$ is a geodesic and $\gamma$ an element of $\Gamma$ such
that $\gamma.c(t) = c(t+1)$ for all
$t \in \Bbb R$; in the terminology of K. Grove \cite{6}, $c$ is called a
$\gamma$-invariant geodesic.
Another such pair
$(c',\gamma')$ represents the same closed geodesic on $Q$ if and only if
there is an element $\overline \gamma \in
\Gamma$ such that $c' = \overline \gamma.c$ and $\gamma' =
\overline \gamma \gamma\ \overline \gamma^{-1}$.

As an example, consider the orbifold $Q$ which is the quotient of the round
2-sphere
$\Bbb S^2$ by a rotation $\rho$ of angle $\pi$ fixing the north pole $N$
and the south pole $S$. The
quotient space $|Q|$ looks like a rugby ball, with two conical points $[N]$
and $[S]$, images of
$N$ and $S$. There are two homotopy classes of free loops on $Q$. Closed
geodesics homotopic to a
constant loop are represented by a closed geodesic on $S^2$ (their length
is an integral multiple
of $2\pi$). If they have positive length, their image in $|Q|$ is either
the equator, a figure
eight or a meridian (image of a great circle through $N$ and $S$). Closed
geodesics in the other homotopy
class are represented by a pair
$(c,\rho)$, where $c$ is either the constant map to $N$ or $S$, or maps
$[0,1]$ to a geodesic arc on the
equator of length an integral odd multiple of $\pi$.

\proclaim{2.4.6. Proposition} Let $Q$ be a compact orbifold.  There is an
$a >0$ such that any
closed geodesic of energy
$<
a$ has zero length. \endproclaim

\demo{Proof} Using the compactness of $|Q|$, one can find a
finite number of convex geodesic balls $X_1,\dots,X_r$ in $X$ such that
 the pseudogroup of change of charts restricted to
$X_i$ is generated by  a finite group
$\Gamma_i$ of isometries of $X_i$, and such that the  union of the $q(X_i)$
is $|Q|$. In particular there
are no closed geodesics of positive length contained in $X_i$.  Using again
the compactness of $|Q|$ one can
find a positive number
$\rho$ such that, for every point $z \in |Q|$, there is a point $x$ in some
$X_i$ with $z =q(x)$ which
is the center of a geodesic ball of radius $\rho$ contained in $X_i$.

 Choose $a < \rho^2/2$. Any closed geodesic
with energy $<a$ is represented by a pair $(c,\gamma)$, where $c: [0,1] \to
X$ is a geodesic
segment of length $< \rho$ contained in some $X_i$ and $\gamma \in
\Gamma_i$ with $\dot c(1) =
\gamma.\dot c(0)$.  If $m$ is the order of $\gamma$, then the $m$-th
iterate $c^m$ of this closed $\Cal
G$-geodesic as defined in 2.3.3, gives a closed geodesic contained in
$X_i$, hence is of
length zero. Therefore $c$ is a comstant map. $\square$
\enddemo

\subheading{2.5 Classifying spaces}

For any topological groupoid $(\Cal G,X)$, one
can construct a classifying space $B\Cal G$, base space of a principal
$\Cal G$-bundle $E\Cal G
\to B\Cal G$. One possible construction  is the geometric realization of
the nerve of the
topological category $\Cal G$. This construction is functorial with respect
to continuous
homomorphisms of groupoids.
 When
$\Cal G$ is the groupoid of germs of changes of charts of an atlas of
uniformizing charts for a
Riemannian  orbifold
$Q$ of dimension
$n$, there is an explicit construction of $B\Cal G$ which is independent of
the particular atlas
defining
$Q$ and which will be therefore noted $BQ$ (see \cite{9} ).

\mn{\bf 2.5.1. Construction of the classifying space $BQ$.}
 Consider
the bundle of
orthonormal coframes $FX$ on $X$; an element of $FX$ above $x \in X$ can be
identified to a
linear isometry from the tangent space $T_xX$ at $x$ to the Euclidean space
$\Bbb R^n $; the
group
$O(n)$ of isometries of $\Bbb R^n$ acts naturally on the left on $FX$ and
this action commutes
with the right action of   the groupoid
$\Cal G$  on $FX$ through the composition with the differential of the
elements of $\Cal G$. As
the action of $\Cal G$ on $FX$ is free, the quotient $FX/\Cal G$ is a
smooth manifold $FQ$
depending only on
$Q$ and not of a particular atlas defining $Q$. The left action of $O(n)$
on $FX$ gives a
locally free action of $O(n)$ on $FX/ \Cal G = FQ$.

 Choose a principal universal $O(n)$-bundle
$EO(n)
\to BO(n)$ for the orthogonal group $O(n)$
 and take for $E\Cal G$ the associated bundle $EO(n) \times_{O(n)}  FX$,
quotient  of $EO(n)
\times FX$ by the diagonal action of $O(n)$. The projection from $FX$ to
$X$ gives a projection
$\alpha_{E\Cal G}: E\Cal G \to X$ with contractible fibers isomorphic to
$EO(n)$. As
the  action of
$O(n)$ on
$FX$  commutes with the natural right action of
$\Cal G$, we get   a free  action of $\Cal G$ on $E\Cal G$ with respect to
the projection
$\alpha_{E\Cal G}$, and
$E\Cal G \to  E\Cal G /\Cal G = EO(n) \times_{O(n)} FX/\Cal G $ is a
universal principal
$\Cal G$-bundle whose base space $B\Cal G = E\Cal G/\Cal G = EO(n)
\times_{O(n)} FQ$ will be
noted
$BQ$. There is a canonical map $\pi:BQ \to |Q|$ induced from the map $q
\circ \alpha_{E\Cal G} :
E\Cal G \to X/\Cal G = |Q|$; it is the projection of the morphism from
$BQ$ to
$\Cal G$ associated to the principal $\Cal G$-bundle $E\Cal G \to BQ$; the
fiber of $\pi$ above a point
$z = q(x)$ is an acyclic space with fundamental group isomorphic to the
isotropy subgroup $\Cal
G_x$ of $x$.

\thickspace

\mn{\bf 2.5.2. What does it classify?} The "maps" from a polyhedron $K$ to
$Q$ (i.e. the morphisms
from
$K$ to
$\Cal G$. see 2.2.1, or equivalently the isomorphisms classes of principal
$\Cal G$-bundle over K,
see 2.2.2)  correspond bijectively (see \cite{9}) to the equivalence
classes of continuous maps
from
$K$ to
$BQ$, two such maps being equivalent if they are connected by an homotopy
along  the fibers of
the projection
 $BQ
\to |Q|$. If $L$ is a subpolyhedron of $K$ and $F$ a principal $\Cal
G$-bundle over $L$,
isomorphic to the pull back of $E\Cal G \to B\Cal G $ by a continuous maps
$f_L: L  \to B\Cal G$, the
morphisms from $K$ to $\Cal G$ relative to $F$ correspond bijectively to
the equivalence classes of
continuous maps $f: K
\to B\Cal G$ whose restriction to $L$ is $f_L$, two such maps being
equivalent if they are connected by an
homotopy along the fibers of the projection to $X/\Cal G = |Q|$ which is
fixed on
$L$ (cf
\cite{9}). This property, called universal property, characterizes a
classifying space for $\Cal G$
up to weak homotopy equivalence.

\thickspace

\mn{\bf 2.5.3. Homology properties.} The projection $BQ \to |Q|$ induces an
isomorphism on rational
homology (or cohomology), because the fibers have trivial rational
homology. If $Q$ is a connected
compact orientable orbifold of dimension $n$, then an orientation
determines a fundamental integral
class which is a generator of
$H_n(|Q|,\Bbb Z)= \Bbb Z$. By the isomorphism $H_n(BQ, \Bbb Q) \cong
H_n(|Q|,\Bbb Q)$, this class
corresponds to a generator of $H_n(BQ,\Bbb Q) =\Bbb Q$ called the
fundamental class of the
oriented orbifold  $Q$.

The projection $BQ \to |Q|$ induces a surjective homomorphism on the
fundamental groups. In
general the homotopy groups of $BQ$, which are isomorphic to the homotopy
groups of $Q$, in the
orbifold sense (see 2.2.3 and 2.5.2), are quite different from the homotopy
groups of $|Q|$.

\heading 3. The free loop space of a Riemannian orbifold \endheading

We consider in this section a Riemannian orbifold $Q$ defined (see
2.1.4) as the quotient
$\Cal G\backslash X$ ($X$ can be the disjoint union of the sources of an
atlas of uniformizing
charts and $\Cal G$ is the groupoid of germs of change of charts). We
recall that
$\Omega^c_X = \bigcup_{x \in X} \Omega^c_x$ is the union of the sets
$\Omega^c_x$
of continuous $\Cal
G$-loops based at $x$. The groupoid $\Cal G$ acts on $\Omega^c_X$ and the
quotient by this action
is the set $|\Lambda^c Q|$ of continuous free loops on $Q$.

\subheading{3.1 The Banach orbifold $\Lambda^c Q$}

\proclaim{3.1.1. Proposition} The set $\Omega^c_X$ of continuus based $\Cal
G$-loops, as well as the set
$\Omega^c_{x,y}$ of equivalence classes of continuous $\Cal G$-paths from
$x$ to $y$,
has a natural structure
of Banach manifold. \endproclaim

\demo{Proof} On $\Omega^c_X$ the structure of Banach manifold is constructed as
follows. Let $c=(g_0,c_1,g_1,\dots,c_k,g_k)$ be a closed $\Cal G$-path over the
subdivision
$0=t_0<t_1<\dots<t_k=1$ based at $x$ (the target of
$g_0$). Let $c^*TX$ be the vector bundle over $\Bbb S^1 = \Bbb R/\Bbb Z$
which is the quotient
of the disjoint union of the bundles $c^*_iTX$ by the equivalence relation
which identifies the
point
$(t_i,\xi_i)
\in c^*_iTX$ with the point $(t_i, Dg_i(\xi_i)) \in c^*_{i-1}TX$ for
$1<i\le k$ and $(0,\xi_0)
\in c^*_1TX$ with $(1,D(g_kg_0)(\xi_0)) \in c^*_kTX$. The projection to the
base space $\Bbb S^1
$ maps the equivalence class of $(t,\xi)$ to $t$ modulo $1$. On the fibers  we
have a scalar product induced from the scalar product given on the fibers
of  $TX$ by the
Riemannian structure on $X$. If $c'$ is a closed $\Cal G$-path based at $x$
equivalent to $c$,
there is a natural isomorphism between $c^*TX$ and ${c'}^*TX$. The tangent
space of $\Omega^c_X$
at
$[c]_x$ will be the Banach space  $C^0(\Bbb S^1,c^*TX)$ of continuous
sections of the
bundle $c^*TX$ with the sup norm. Such a section $v$ is represented by a
sequence
$(v_1,\dots,v_k)$, where $v_i$ is a vector field along $c_i$ such that the
above compatibility
conditions  are satisfied  at the points $t_i$.

	Given $c$, choose $\epsilon >0$ so small that,  for each $t \in
[t_{i-1},t_i]$, the
exponential map $\exp_{c_i(t)}$ is defined and is injective on the ball of
radius $\epsilon$
in $T_{c_i(t)}X$; then $g_i$ extends uniquely (see 2.1.5) to a section
$\tilde g_i$ of
$\alpha$ defined on the ball of radius $\epsilon$ and center $\alpha
(g_i)$. Let
$\tilde U_c^\epsilon$ be the open ball of radius $\epsilon$ in $C^0(\Bbb
S^1,c^*TX)$ centered at the origin. We define
the chart $\exp_c^\epsilon :\tilde U_c^\epsilon
\to U_c^\epsilon \subseteq \Omega_X$  by mapping a section $v$ given by
$(v_1,\dots,v_k)$ to the
equivalence class of the based closed $\Cal G$-path $c^v=
(g^v_0,c^v_1,g^v_1, \dots,
c^v_k,g^v_k)$ defined as follows: $c^v_i(t) = \exp_{c_i(t)}v_i(t)$ ,
$g^v_{i} = \tilde
g_{i}(c^v_{i+1}(t_i)$ for $i<k$ and $g^v_k=\tilde g_k(\omega(g^v_0))$.

 It is easy to see that each $\exp_c^\epsilon$ is a bijection (see 2.3.2),
that the
images
$U_c^\epsilon$ for various $c$ and $\epsilon$ form a basis for a topology
on $\Omega^c_X$, and
that the change of charts are differentiable.

Therefore
$\Omega^c_X$ is a Banach manifold.
The
tangent space $T_{[c]_x}\Omega^c_X$ at $[c]_x$ is the space of continuous
sections of the
vector bundle
$c^*TX$ over $\Bbb S^1$. Note that $\Omega^0_X$ is a finite dimensional
submanifold of
$\Omega^c_X$.

The Banach manifold structure on $\Omega^c_{x,y}$ is defined similarly. For a
$\Cal G$-path
$c=(g_0,c_1,\dots, c_k,g_k)$, the tangent space at $[c]_{x,y}$ is
isomorphic to the space of
vector fields
$v= (v_1,\dots,v_k)$ along $c$ which vanish at $0$ and $1$. $\square$ \enddemo

The action of $\Cal G$ on
$\Omega^c_X$ with respect to the projection assigning to  a based $\Cal
G$-loop its
base point is differentiable. The quotient of $\Omega^c_X$ by this action is by
definition the "space"
of (continuous) free loops $|\Lambda^c(\Cal G)| = |\Lambda^c Q|$  on $Q$.

\proclaim{3.1.2. Proposition} The space $|\Lambda^c Q|$ of continuous free
$\Cal
G$-loops on $Q$ has a natural Banach orbifold
structure noted
$\Lambda^c Q$. The subspace $|\Lambda^0 Q|$ of free loops of length zero is a
"suborbifold"
$\Lambda^0 Q$ of $\Lambda^c Q$. \endproclaim

\demo{Proof} If $q_i: X_i \to V_i$ is a uniformizing chart for $Q$ (see
2.1.1), then
$\Omega^c_{X_i} \to
\Gamma_i \backslash \Omega^c_{X_i}$ is a uniformizing chart for $\Lambda^c Q$,
where $\Omega^c_{X_i} =
\bigcup_{x \in X_i} \Omega^c_x$. The groupoid of germs of changes of chart is
the groupoid
$\overline{\Cal G}:=\Cal G
\times_X
\Omega^c_X$, the subspace of
$\Cal G \times
\Omega^c_X$ consisting of pairs $(g,[c]_x)$ with $\alpha(g) = x$. The source
(resp. target)
projection maps $(g,[c]_x)$ to $[c]_x$ (resp. $^g[c]_x$). The composition
$(g',[c']_{x'})(g,[c]_x)$ is defined if  $^g[c]_x =[c']_{x'}$ and is equal
to $(g'g,[c]_x)$.

The suborbifold structure on $\Lambda^0 Q$ is obtained by replacing
$\Omega^c_X$ by $\Omega^0_X$ and
$\Cal G
\times_X
\Omega^c_X$ by $\Cal G \times_X \Omega^0_X$. Note that $\Lambda^0Q$ was
considered by
Kawasaki in \cite{10}. $\square$
\enddemo

\mn{\bf 3.1.3. The action of $\Bbb S^1$ on the orbifold $\Lambda^c Q$.} The
action
of $\Bbb S^1$ on $|\Lambda^c Q|$ described at the end of 2.3.4 comes from a
continuous action of $\Bbb S^1$ on the orbifold $\Lambda^c Q$. This means
that if
$c$ is a continuous $\Cal G$-loop based at $x$ and if $c'$ is a  $\Cal
G$-loop based at $x'$ representing the translate $\tau.[c]$ of $[c]$ by
$\tau \in
\Bbb R
\backslash \Bbb Z = \Bbb S^1$, there is a diffeomorphism $h_\tau$ of a
neighbourhood of $[c]_x$ to a neighbourhood of $[c]_{x'}$ projecting to the
translation by
$\tau$ of a neighbourhood of $[c]$ to a neighbourhood of $[c']$ depending
continuously on $\tau$. As explained for instance in \cite{12, p. 39}, in
the charts
$\exp^\epsilon_c$ and $\exp^\epsilon_{c'}$, indeed  $h_\tau$ is given by a
linear
map. The continuity in $\tau$ is checked as in \cite{12}  using a chart
associated
to a $C^\infty$ based $\Cal G$-loop close to $[c]_x$.

\subheading{3.2 A classifying space for $\Lambda^c Q$}

We consider as in 2.5 a classifying space $BQ$, base space of a universal
principal $\Cal
G$-bundle
$E\Cal G \to B\Cal G = BQ$. Let $ E\Cal G \times_X \Omega^c_X$ be the
subspace of $ E\Cal G \times
\Omega^c_X$ consisting of pairs $(e,[c]_x)$ such that $\alpha_{E\Cal G}(e)
= x$.
We note $E\Cal G
\times_\Cal G
\Omega^c_X$ its quotient by the equivalence relation identifying
$(e.g,[c]_x)$ to $(e, ^g[c]_x)$. We can consider $ E\Cal G
\times_X\Omega^c_X \to
E\Cal G \times_\Cal G \Omega^c_X$ as a
principal
$(\Cal G \times_X \Omega^c_X)$-bundle, the action map being the natural
projection to
$\Omega^c_X$.

\proclaim{ 3.2.1.  Proposition} $E\Cal G \times_X \Omega^c_X \to E\Cal G
\times_\Cal G \Omega^c_X$ is a
principal universal
$(\Cal G \times_X \Omega^c_X)$-bundle. The base space $E\Cal G \times_\Cal G
\Omega^c_X$
will be  noted $B\Lambda^c Q$.

 Similarly $E\Cal G \times_X \Omega^0_X \to E\Cal G \times_\Cal G
\Omega^0_X$ is a
principal universal
$(\Cal G \times_X \Omega^0_X)$-bundle. Its base space $E\Cal G \times_\Cal
G \Omega^0_X$
is  noted $B\Lambda^0 Q$.

	The natural projection $B\Lambda^c Q \to BQ$ (i.e. $E\Cal G
\times_\Cal G \Omega^c_X \to E\Cal
G/\Cal G = B\Cal G = BQ$) is a Serre fibration with fibers isomorphic to
$\Omega^c_x$.

\endproclaim

\demo {Proof} With respect to the action map to $X$,  $E\Cal G$ is
a locally
trivial bundle  with contractible fibers. The pull back of this bundle by
the projection
$\Omega^c_X \to X$ is the bundle
$E\Cal G
\times_X \Omega^c_X$ with base space $\Omega^c_X$. It has also contractible
fibers and therefore
$E\Cal G \times_X \Omega^c_X \to E\Cal G \times_\Cal G \Omega^c_X$ is a
universal principal $(\Cal G
\times_X
\Omega^c_X)$-bundle. The same argument works with $\Omega^0_X$ replacing
$\Omega^c_X$. It is easy to see
that the projection
$\Omega^c_X
\to X$ is a Serre fibration. Therefore the projection
$ E\Cal G
\times_\Cal G
\Omega^c_X \to B\Cal G$ is also a Serre fibration because it is locally a
pull back of the
fibration
$\Omega^c_X \to X$, as it is seen using a local section of the projection
$E\Cal G \to B\Cal
G$. $\square$
\enddemo

Let $\Lambda BQ$ be the space of continuous free loops on $BQ$.

\proclaim{3.2.2. Theorem} $\Lambda BQ$  is the base space of a universal $(\Cal
G \times_X \Omega^c_X)$-bundle. This bundle is the pull back of $E\Cal G
\times_X \Omega^c_X$ by a map
$\phi: \Lambda BQ \to B\Lambda^c Q$ which is a weak homotopy equivalence and
commutes with the
projections to $BQ$. Therefore $\Lambda BQ$ is a classifying space for the
orbifold $\Lambda^c Q$.

For points $z \in BQ$ and $x \in X$ projecting to the same point of $|Q|$,
the map $\phi$ induces
a weak homotopy equivalence from the space $\Omega_z BQ$ of loops on $BQ$
based at $z$ to
$\Omega^c_x$.\endproclaim

We first prove two lemmas. The first one is a tautology.

\proclaim{3.2.3. Lemma} There is a canonical principal $\Cal G$-bundle
$E_{\Omega^c_X}$ over $\Omega^c_X
\times
\Bbb S^1$ relative to the bundle over $\Omega^c_X \times \{1\}$ pull back of
the bundle $\Cal G \to X$
by the projection $([c]_x, 1) \mapsto x$ to $X$. Its restriction to
$\{[c]_x\} \times \Bbb S^1$ is
the relative principal $\Cal G$-bundle $E_{[c]_x}$  over $\Bbb S^1$
associated to $[c]_x$ (cf. 2.3.2).
\endproclaim

\demo{Proof} Locally the bundle is constructed above the product of a
neighbourhood $U$ of $[c]_x$
with $\Bbb S^1$ as in 2.3.2 using a chart, and those local constructions
are uniquely glued together
uisng 2.3.3. $\square$\enddemo

\proclaim{3.2.4. Lemma} Let $K$ be a topological  space and $L$ be a
subspace of $K$. Let $f_0: L
\to
\Omega_x$ be a continuous map and let $F$ be the principal $\Cal G$-bundle
over
 $A:=( K
\times \{1\}) \cup  (L \times \Bbb S^1) \subseteq K \times \Bbb S^1$ whose
restriction to $L
\times \Bbb S^1$ is the  pull back of
$E_{\Omega^c_X}$ by the map
 $f_0 \times id$ and whose restriction to $K \times \{1\}$ is the pull back
from $\Cal G$ by
the constant map to $x$.

There is a bijective correspondence between continuous maps $f: K \to
\Omega^c_x$ extending
$f_0$ and isomorphism classes of principal $\Cal G$-bundle over $K \times
\Bbb S^1$ relative to $F$.
\endproclaim

\demo{Proof} Given $f$, it is clear that the pull back of $E_{\Omega^c_x}$ by
$f \times id$ is a
 bundle relative to $F_0$. Conversely, given such a bundle $E$, for each $y
\in K$, its
restriction to
$\{y\} \times \Bbb S^1$ gives a $\Cal G$-loop based at $x$ that we note
$f(y)$. It is clear that
the map $f: K \to \Omega^c_x$ extends $f_0$. $\square$ \enddemo

\demo{Proof of the theorem} Let $\Cal U = \{U_i\}_{i\in I}$ be an open
cover of $B\Cal G = BQ$
with local
 sections $s_i:U_i \to E\Cal G$ with respect to the projection $E\Cal G \to
B\Cal G$. Let
$f_{ij}:U_i \cap U_j \to \Cal G$ be the corresponding cocycle: for $z \in
U_i \cap U_j$, we
have $s_j(z) = s_i(z)f_{ij}(z)$; in particular, $f_i:= f_{ii}: U_i \to X$
is equal to
$\alpha_{E\Cal G}
\circ s_i$.

	Let $\pi: \Lambda B\Cal G \to B\Cal G$ be the map associating to a
loop $l: [0,1] \to B\Cal G$
its base point $l(0)=l(1)$. For each
$i \in I$, we consider the continuous map $\phi_i: \pi^{-1}(U_i) \to
\Omega_X$ defined as
 follows. Given a loop $l$ on $BQ$ based at $z \in U_i$, we choose a
subdivision $0=t_0
<t_1<\dots<t_k=1$ of
$[0,1]$ such that
$l([t_{j-1},t_j])
\subseteq U_{r_j}$ for some $r_j \in I$. Then $\phi_i(l)$ is  the $\Cal
G$-loop based at $f_i(z)$ represented by the $\Cal G$-path
$c=(g_0,c_1,\dots,c_k,g_k)$ where
$c_j:[t_{j-1},t_j] \to X$ maps
$t$ to $f_{r_j}(l(t))$, $g_j = f_{r_jr_{j+1}}(l(t_j))$ for $1 \le j \le k-1$,
$g_0=f_{ir_1}(l(0))$ and $g_k=f_{r_ki}(l(1))$. It is easy to see that the
equivalence class of
$c$ does not depend of the choices and that the map $\phi_i$ is continuous.

By construction the restriction of the map $\phi$ to $\pi^{-1}(U_i)$ is
defined as the
composition of $(s_i,\phi_i): \pi^{-1}(U_i) \to E\Cal G \times_X \Omega^c_X $
with the
projection to the quotient $E\Cal G \times_\Cal G \Omega^c_X = B\Lambda^c Q$.

	 The spaces
$\Lambda BQ$ and $B\Lambda^c
Q$ are Serre fibrations over $BQ$ and $\phi$ commutes with the projections.
To prove that $\phi$ is a weak homotopy equivalence it is sufficient to
check that, for a base point
$z
\in U_i$, the restriction of $\phi$ to the space $\pi^{-1}(z) =
\Omega_z(BQ)$ of loops on $BQ$
based at
$z$ is a weak homotopy equivalence to the fiber of $E\Cal G \times_\Cal G
\Omega^c_X \to
BQ$ above $z$. This fiber is isomorphic to $\Omega^c_x$, the space of
continuous
$\Cal G$-loops based at
$x=f_i(z)$, by the map sending $[c]_x$ to the class modulo $ \Cal G$ of
$(s_i(z),[c]_x)$. Via this
isomorphism, $\phi$ maps a loop $l \in \Omega^c_z$ to $\phi_i(l)$.

Fix a base point $l \in \Omega_z BQ$ and let $[c]_x$ be its image by
$\phi_i$. We shall apply below
the lemma 3.2.4 with
$K=\Bbb S^m$,
$L=*$ a base point and
$f_0$ mapping
$*$ to
$[c]_x$. We first observe that the set $\pi_m(\Omega_zBQ;l)$ is in
bijection with the set of
homotopy classes of maps from $\Bbb S^m \times \Bbb S^1$ to $BQ$ which maps
$\Bbb S^m \times
\{0\}$ to $z$ and whose restriction to $\{*\}
\times \Bbb S^1$ is equal to $l$. According to 2.5.2, this set is in
bijection with the set of
homotopy classes of principal $\Cal G$-bundles over
$\Bbb S^m \times \Bbb S^1$ relative to
$F$. Using 3.2.4, we see that this set corresponds bijectively  to the set
of homotopy classes
of maps
$(\Bbb S^m,*)
\to (\Omega^c_x,[c]_x)$.
$\square$ \enddemo

\remark{3.2.5. Remark}
Composing $\phi: \Lambda BQ \to B\Lambda^c Q$ with the natural projection from
$B\Lambda^c Q$ to
$|\Lambda^c
Q|$, we get an $\Bbb S^1$-equivariant map
$$\Lambda BQ \to |\Lambda^c Q|$$ with respect to the natural action of
$\Bbb S^1$ on
free loops.\endremark
\remark{3.2.6. Remark} One can show that, for a space $K$,  there is a
canonical correspondance associating
to a principal $\Cal G$-bundle $E$ over $K \times \Bbb S^1$ a principal
$\overline{\Cal G}-bundle$ over $K$,
where $\overline{\Cal G} = \Cal G \times_X \Omega^c_X$, inducing a bijection
on isomorphisms classes. The
universal $\overline{\Cal G}$-bundle over $\Lambda BQ$ corresponds to the
principal $\Cal G$-bundle over
$\Lambda BQ \times \Bbb S^1$ which is the pull back of $E\Cal G$ by the
evaluation map $\Lambda BQ \times
\Bbb S^1 \to BQ$ sending $(l,t)$ to $l(t)$. \endremark

\subheading{3.3 The Riemannian  orbifold $\Lambda Q$ of free $\Cal G$-loops
of class $H^1$}

We consider as above a Riemannian orbifold $Q = \Cal G \backslash X$.

\mn{\bf 3.3.1. $\Cal G$-paths of class $H^1$.} A $\Cal G$-path $c=
(g_0,c_1,g_1,\dots, c_k,g_k)$  over a subdivision $0=t_0<t_1<\dots<t_k=1$
is of class
$H^1$ if each $c_i$ is absolutely continuous and the velocity functions $t
\mapsto |\dot c_i(t)|$ are
square integrable. Those conditions are also satisfied by any $\Cal G$-path
in the equivalence
class of
$c$. We denote $\Omega_{x,y}$ the set of equivalence classes of $\Cal
G$-paths of class $H^1$ from
$x$ to
$y$ and $|\Lambda Q|$ (resp.
$\Omega_x$ ) the set of free
$\Cal G$-loops on
$Q$ (resp. based at
$x$) represented by closed $\Cal G$-path of class
$H^1$. Let $\Omega_X= \bigcup_{x \in X}\Omega_x$. The energy function $E$
is defined on all those
spaces.

\medspace

\mn{\bf 3.3.2. $\Omega_{x,y}$ and $\Omega_X$ as Riemannian Hilbert
manifolds.}

On $\Omega_X$, we define a structure of Hilbert Riemannian manifold
using the following charts. Let $c=(g_0,c_1,g_1,\dots,c_k,g_k)$ be a closed
$\Cal
G$-path over the subdivision $0=t_0 < t_1<\dots<t_k=1$ ; we assume that $c$
is piecewise
differentiable, i.e. that each $c_i$ is differentiable. We consider as in
3.1.1 the vector bundle
$c^*TX$ over $\Bbb S^1$ which is the union of the vector bundles $c_i^*TX$
on which we consider
the connection induced from the Levi Civita connection on $TX$. A
continuous section
$v=(v_1,\dots,v_k)$ of $c^*TM$ is an $H^1$-section if each $v_i$ is
absolutely continuous and
the covariant derivative $\nabla v_i$ is square integrable. The space of
$H^1$-sections is a
separable  Hilbert space noted $H^1(c^*TX)$: if $w=(w_1,\dots,w_k)$ is
another $H^1$-section, then
the scalar product is defined by the following formula
$$(v,w) = \sum_i \int_ {t_{i-1}}^{t_i}[\langle v_i(t),w_i(t)\rangle +
\langle\nabla v_i(t),\nabla
w_i(t)\rangle]dt.\tag{3.3.3}$$ This definition is independent of the
particular choice of $c$ in its equivalence class.

Consider the open set of $H^1(c^*TM)$ consisting of sections $v$ such that
$|v_i(t)| < \epsilon$, where $\epsilon$ is like in 3.1.1. The map $\exp_c$
maps this open set to a
subset ${U'}_c^\epsilon$ of $\Omega_X$. When $c$ varies over the  closed
piecewise
differentiable $\Cal G$-paths, one proves like in the classical case of
Riemannian manifolds (see
for instance \cite{11,12  }) that one obtains in this way  an atlas defining on
$\Omega_X$ a structure of  Hilbert manifold. When $c$ is a closed $\Cal
G$-path based at $x$ which
is of class
$H^1$, then the bundle $c^*TX$ has also a natural connection induced
from the connection on $TX$, and one can also define, using the formula
3.3.3,  the Hilbert space
$H^1(c^*TX)$ of
$H^1$-sections of this vector bundle (namely sections $v=(v_1,\dots,v_k)$
such that each $v_i$
gives a map from $[t_{i-1},t_i]$ to $TX$ which is of class $H^1$); it is
canonically isometric to
the tangent space of
$\Omega_X$ at the equivalence class $[c]_x$ of $c$.

Similarly the set $\Omega_{x,y}$  is
naturally a Riemannian manifold, the tangent space at $[c]_{x,y}$ is the
space of $H^1$-sections
of the bundle $c^*(TX)$ over $[0,1]$ which vanish at $0$ and $1$, with the
scalar product given by
3.3.3. In particular
$\Omega_x$ is a closed submanifold of $\Omega_X$. It is complete as a
Riemannian manifold if the
orbifold $Q$ is complete.

\medskip\noindent{\bf 3.3.4. The Riemannian orbifold $\Lambda Q$.} Like in
3.1.2 the
space $|\Lambda Q|=
\Cal G \backslash \Omega_X$ of free $\Cal G$-loops on $Q$ of class $H^1$ as
a natural
Riemannian orbifold structure noted $\Lambda Q$. The groupoid of germs of
change of charts for
$\Lambda Q$ is
$\Cal G
\times_X
\Omega_X$. As in 3.2.1 one proves that the  space $B\Lambda Q :=E\Cal G
\times_\Cal G \Omega_X$
is a classifying space for
$\Lambda Q$.

 Like in 3.1.3  one can check (see \cite{12, p.39} that the action of the group
$\Bbb S^1$ on $|\Lambda Q|$ comes from an action of $\Bbb S^1$ on the orbifold
$\Lambda Q$  by isometries which is continuous (but not differentiable).

\proclaim {3.3.5 Proposition} The natural inclusions
$$\Omega_X \to \Omega^c_X,\ \ \ \  \Omega_{x,y} \to \Omega^c_{x,y}$$  are
continuous and are
homotopy equivalences. In particular, if $Q$ is connected, $\Omega_{x,y}$
has the same weak homotopy
type as the space of loops on $BQ$ based at a fixed point.

The induced inclusion
$$B\Lambda Q \to B\Lambda^c Q$$ is a  homotopy equivalence. \endproclaim

\demo{Proof}  Denote by $P$ the
spaces
$\Omega_X$ or
$\Omega_{x,y}$ and by
$P^c$ the spaces $\Omega^c_X$ or $\Omega^c_{x,y}$. The fact that the
inclusions $i:P \to P^c$ are continuous is
proved like in the classical case. To show
 that they are homotopy equivalences, we follow the argument of Milnor
\cite{13, p. 93-94}.  For a positive
integer $k$, let $P_k$ (resp.$P^c_k$) be the subspace of
$P$ (resp.
$P^c$) formed by the elements represented by $\Cal G$-paths
 $c=(g_0,c_1,g_1,\dots, c_{2^k},g_{2^k})$ defined over the subdivision $0=t_0
<\dots<t_{2^k}=1$, where $t_i =i/2^k$, and such each $c_i(t_i)$ is the
center of a convex geodesic ball
containing the image of $c_i$. One can deform continuously (see \cite{13,
p. 91}) such a $\Cal
G$-path $c$ to the $\Cal G$-path $\overline c =(g_0,\overline
c_1,g_1,\dots, \overline c_{2^k}, g_{2^k})$,
where
$\overline c_i$ is the geodesic segment joining $c_i(t_{i-1})$ to
$c_i(t_i)$.  Passing to equivalence
classes, this gives a continuous deformation of $i_{|P_k}: P_k \to P^c_k$
and implies that this inclusion is
a homotopy equivalence. As the spaces $P$ and $P^c$ are the increasing union
of the open
subspaces $P_k$ and $P^c_k$ for $k=1,2,\dots$, it follows that $i$ is a
homotopy equivalence (see
the appendix of \cite{13}).

The last assertion follows from the fact that the above deformation
commutes with the projection to $X$ and
with the action of $\Cal G$.
$\square$
\enddemo

\mn{\bf 3.3.6. The developable case.} Let $X$ be a simply connected
Riemannian manifold,
$\Gamma$ a discrete subgroup of the group of isometries of $X$, and let $Q$
be the quotient
orbifold $\Gamma \backslash X$.  For $\gamma \in \Gamma$, the space of
$H^1$-curves $c:[0,1]
\to X$ such that $c(1) =\gamma.c(0)$ is noted $\Lambda(X,\gamma)$ in
Grove-Tanaka \cite{6}. It
is connected because $X$ is assumed to be simply connected. The space
$\Omega_X=
\Omega_X(\Gamma \ltimes X)$ of based $\Gamma \ltimes X$-loops of class
$H^1$ is the disjoint
union  $$\Omega_X =\coprod_{\gamma \in  \Gamma}\Lambda(X,\gamma).$$
The connected components of the space $|\Lambda Q|$ of free loops of class
$H^1$ are in
bijection with the conjugacy classes in $\Gamma$. To the conjugacy class of
$\gamma\in \Gamma$
corresponds the quotient  of $\Lambda(X,\gamma)$ by the action of the
centralizer $Z_\gamma$
of $\gamma$ in $\Gamma$.

Grove and Tanaka consider the case where $X$ is compact simply connected.
In that case
$\Gamma$ is a finite group. For $\gamma \in \Gamma$, the existence of one
(resp. infinitely
many geometrically distinct ) $\gamma$-invariant geodesic of positive
length is equivalent to
the existence of one (resp. infinitely many geometrically distinct) closed
geodesic on $Q$ of
positive length in the connected component $Z_\gamma \backslash \Lambda
(X,\gamma)$ of
$\Lambda Q$.

\mn{\bf 3.3.7  Classification of tubular neighbourhoods of $\Bbb
S^1$-orbits.} We
consider an effective continuous action of
$\Bbb S^1$ by isometries on a Riemannian orbifold. We have in mind the case
of the natural action by isometries of $\Bbb S^1$ on $\Lambda Q$. We want
to describe an invariant  tubular neighbourhood of a smooth
$\Bbb S^1$-orbit of positive length. We can assume that the orbifold is
developable because such a
tube is always developable. So we consider an orbifold $Q$ which is the
quotient of a
connected Riemannian manifold $Y$ by a discrete group $\Gamma$ of
isometries of $Y$ acting properly. We
assume that we have a continuous effective action of $\Bbb S^1$ on $Q$ by
isometries. This is
equivalent to say that we have a continuous action of $\Bbb R$ on $Y$ by
isometries which commutes with
the action of
$\Gamma$ and such that there is a unique element $\gamma_0 \in \Gamma$ such
that, for each
$y \in Y$, the translate $T_1y$ of $y$ by $1 \in \Bbb R$ is equal to
$\gamma_0.y$. This implies that $\gamma_0$ is in the center of $\Gamma$. We
note below
$T_\tau y$ the translate of
$y$ by
$\tau
\in
\Bbb R$.

Let $y_0$ be a point of $Y$ which is not fixed by the action of $\Bbb R$,
and let $\gamma_0$ be as
above. We assume that the $\Bbb R$-orbit of $y_0$ is differentiable. Let
$\overline \Gamma$ be the
subgroup of $\Gamma$ leaving invariant the $\Bbb R$-orbit of $y_0$.
Consider the subgroup
$G$ of
$\overline \Gamma \times \Bbb R$ consisting of pairs $(\gamma,\tau)$ such
that $\gamma.y_0 = T_\tau y_0$.
The image of $G$ by the natural projection $(\gamma,\tau) \mapsto \tau$ is
a discrete subgroup of $\Bbb R$
generated by a number $1/r$, where $r$ is a positive integer. The kernel of
this projection is
canonically isomorphic to the finite subgroup $\Gamma_0$ of
$\overline\Gamma$ fixing $y_0$. We have an
exact sequence
$$ 1 \to \Gamma_0 \to G \overset\psi\to\longrightarrow \Bbb Z\ 1/r  \to 0.$$
Let $B$ be the image by the exponential map $\exp_{y_0}$ of a small ball
centered at $0$ in the
subspace of the tangent space at $y_0$ orthogonal to the velocity vector of
the $\Bbb R$-orbit of
$y_0$. In particular
$B$ is a slice for the
$\Bbb R$-action and is left invariant by $\Gamma_0$.

\proclaim{3.3.8. Proposition} Let $I$ be the quotient of $G$ by the
subgroup generated by
$(\gamma_0,1)$. The group $I$ acts effectively on $B$ and the homomorphism
$\psi$ gives
(modulo $1$) a homomorphism from
$I$ to the subgroup of order $r$ of $\Bbb S^1$.
A
$\Bbb S^1$-invariant tubular neighbourhood $U$ of the
$\Bbb S^1$-orbit of
$[y_0]$ in $|Q|$ is $\Bbb S^1$-equivariantly isomorphic as an orbifold to
the quotient  $\Bbb S^1 \times_I B$
of $
\Bbb S^1
\times B$ by the diagonal action of $I$. \endproclaim

\demo{Proof}
The group $G$ acts on $B$ (resp. $\Bbb R$), the element $(\gamma,\tau) \in
G$ mapping
$b \in B$ to $\gamma
T_\tau^{-1}(b)$ (resp. $t$ to $t+\tau$). Note that $(\gamma_0,1)$ acts
trivially on $B$.
Consider the immersion
$i:
\Bbb R \times B
\to Y$ mapping
$(t,b)$ to
$T_t b$. Its image is an $\Bbb R$-invariant tubular neighbourhood
$\overline U$ of the $\Bbb R$-orbit of $y$.
The group
$G$ acts on $\Bbb R \times B$ by the diagonal action. The
immersion
$i$ is equivariant with respect to the natural projection $G
\to \overline\Gamma$ and with respect to the action of $\Bbb R$ by translations.The immersion $i$ induces an isomorphism from the orbifold quotient of
$\Bbb R \times B$ by $G$ to the orbifold quotient of $\overline U$ by
$\overline \Gamma$ which is an invariant
tubular neighbourhood of the
$\Bbb S^1$-orbit of the projection $[y_0]$ of $y_0$ to $|Q|$. In particular
$G$ is the fundamental group of
this tubular neighbourhood.

The  quotient of $ \Bbb R \times B$ by the subgroup generated by
$(\gamma_0,1)$ is isomorphic
to
$\Bbb S^1
\times B$ and its quotient  by $G$ is isomorphic to the quotient $U = \Bbb
S^1 \times_I B$ of $\Bbb S^1
\times B$ by the diagonal action of the group $I$ quotient of $G$ by the
subgroup generated by
$(\gamma_0,1)$. This describes an orbifold $\Bbb S^1$-invariant tubular
neighbourhood of the $\Bbb
S^1$-orbit of $[y_0]$. The action of $I$ on $B$ is effective because the
action of $\Bbb S^1$ on $Q$ is
assumed to be effective.
$\square$
\enddemo

\subheading{3.4 The energy function}

The energy function $E$
is well defined on $\Omega_{x,y}$ and
$\Omega_X$. As it is invariant by the action of $\Cal G$, it gives a well
defined function on
$|\Lambda Q|$ still noted $E$. We list below some of its properties,
refering for instance to Klingenberg
\cite{11,12}  for the proofs.

\mn{\bf 3.4.1.} $E$ is a differentiable function on
$\Omega_X$ or $\Omega_{x,y}$. The gradient
$\text{grad } E$ of $E$ at
$[c]_x$ is given by the formula:
$$\text{grad } E(v) = \sum_{i=1}^k \int_{t_{i-1}}^{t_i} \langle\dot
c_i(t),\nabla v_i(t)\rangle dt,$$
where $v=(v_1, \dots,v_k) \in H^1(c^*TX)$.

The critical points of $E$ on $\Omega_X$ (resp. $\Omega_{x,y}$) are in
one to one
correspondence with the equivalence classes of based closed geodesic $\Cal
G$-paths (resp. of geodesic
$\Cal G$-paths from
$x$ to
$y$).

\mn{\bf 3.4.2. The Palais-Smale condition (C).} For a compact orbifold
(resp. a complete
orbifold) the Palais-Smale condition (C) holds for the function $E$ on
$|\Lambda Q|$ (resp. on $\Omega_{x,y}$).

 Namely, let $c_m$ be a sequence of closed $\Cal G$-paths
(resp. of $\Cal G$-paths from $x$ to $y$) of class $H^1$ such that

(i) the sequence $E(c_m)$ is bounded,

(ii) the sequence $|\text {grad }E(c_m)|$ tends to zero;

\noindent then the sequence $[c_m]$ (resp. $[c_m]_{x,y}$) has accumulations
points and any
converging subsequence converges to a geodesic.

This implies that, for $a \ge 0$, the set of critical points of $E$ in the
subspaces
$E^{-1}([0,a])$ of $|\Lambda Q|$ or $\Omega_{x,y}$ is compact.

\mn{\bf 3.4.3. The gradient flow.} The vector field $-\text{grad } E$
generates a local flow
$\phi_t$ on $\Omega_X$ which commutes with the action of $\Cal G$. On the
quotient
$|\Lambda Q| = \Cal G\backslash \Omega_X$ it gives a local flow which is
defined for all $t\ge 0$ if $Q$ is compact.

Similarly the vector field $-\text{grad } E$ on
$\Omega_{x,y}$ generates a local flow $\phi_t$ which is defined for all $t
\ge 0$ when $Q$ is complete.

On $\Omega_X$, when $Q$ is compact, the local flow $\phi_t$ is
$\overline{\Cal G}$-complete in
the following sense,  where $\overline{\Cal G}: = \Cal G \times_X
\Omega_X$ is the groupoid of
germs of change of charts of the orbifold $\Lambda Q$. Given
$\overline x
\in
\Omega_X$ and
$\tau
\ge 0$, one can find   a
$\overline{\Cal G}
$-path $\overline c= (\overline g_0, \overline c_1, \overline g_1,\dots,
\overline
c_k, \overline g_k)$ over a subdivision $0=t_0 \le t_1
\le \dots  \le t_k=\tau$ of the interval
$[0,\tau]$ such that $\omega(\overline g_0) = \overline x$ and   $\overline
c_i(t) =
\phi_{t-t_{i-1}}(\overline c_i(t_{i-1}))$ for $i=1,\dots, k$. The image of
this path in $|\Lambda Q|$ is
the
$\phi_t$ trajectory of the projection of $\overline x$ for $t \in
[0,\tau]$. When $\overline x$
remains in a small neighbourhood, such a $\overline{\Cal G}$-path exists
over the same subdivision
and varies continuously.

 Given another such  $\overline{\Cal G}$-path
$\overline c'= (\overline g'_0, \overline c'_1, \overline g'_1,\dots, \overline
c'_k, \overline g'_k)$ issuing from
$\overline x'$ defined over the same subdivision of the interval $[0,\tau]$
and an element
$\overline g \in \overline{\Cal G}'$ with source $\overline x$ and target
$\overline x'$, then
there are unique continuous maps  $ h_i:[t_{i-1},t_i] \to \overline{\Cal
G}$ such that
$\alpha( h_i(t)) =\overline c_i(t)$, $\omega(h_i(t))= \overline c'_i(t)$,
$\overline g'_0 h_1(0)=
\overline g \overline g_0$ and $h_i(t_i)\overline g_i = \overline
g'_ih_{i+1}(t_i)$ for
$i=1,\dots, k-1$. This implies easily the following

\proclaim{3.4.4. Lemma}: Given a "map"  $f$ from a compact space $K$ to
$\Lambda Q$
(i.e. a morphism from $K$ to
$\overline {\Cal G}$) and
$\tau
\ge 0$, there is a unique  homotopy $f_t$ of $f$ parametrized by  $t \in
[0,\tau]$ whose
projection  to
$|\Lambda Q|$ is the flow $\phi_t$ applied to the projection of $f$.
\endproclaim.

\mn{\bf 3.4.5 Definition of  $\phi$-families.} Let $P$ be either
$|\Lambda Q|$ or $\Omega_{x,y}$. In the
first case we assume that $|Q|$ compact and in the second case that $|Q|$
is complete.  For a number
$a
\in \Bbb R$, we denote $P^a$ the set of point of $P$ for which the value of
the energy function
is $\le a$.

A $\phi$-family (see \cite{12}, p. 20) is a collection $\Cal F$ of
non-empty subsets $F$ of $P$ such
that
$E$ is bounded on each $F$ and
 for $F \in \Cal F$, then $\phi_t(F) \in \Cal F$ for all $t >0$.

Let $a \in \Bbb R$ and assume that there is $\epsilon >0$ such that $E$ has
no critical values in
$]a +\epsilon]$. A $\phi$-family of $P$ mod $P^a$ if a $\phi$-family $\Cal
F$ such that each
member $F$ is not contained in $P^{a+\epsilon}$. Note that a $\phi$-family
is always a
$\phi$-family  of $P$ mod $P^a$ for $a <0$.

The proof of the following theorem is similar to the proof of the
corresponding theorem in the
classical case (see for instance \cite{12}, p. 21). Up to the end of this
section, we assume that
$Q$ is compact (resp. complete) if $P = |\Lambda Q|$ (resp. $P =
\Omega_{x,y}$).

\proclaim{3.4.6. Theorem} The critical value $a_\Cal F$ of a $\phi$-family
$\Cal F$ of $P$ mod
$P^a$, i.e. the number  defined by
$$a_\Cal F = \inf_{F \in \Cal F} \sup E|F,$$ is always $>a$ and there is a
critical point of $E$ with
value $a_\Cal F$. \endproclaim

Applying this theorem to the $\phi$-family formed by the points of a
connected component of $P$
and for $a <0$, we get the following corollary.
\proclaim{3.4.7. Corollary} The energy function $E$ restricted to a
connected component of $P$
assumes its infimum in some point,  and  such a point  is a critical point
of $E$. \endproclaim

\mn{\bf 3.4.8. General remarks.} Let $c$ be a $\Cal G$-geodesic path from
$x$ to $y$ (respectively a closed
$\Cal G$-geodesic path based at $x$). The usual expression of the Hessian
of the energy function at
$T_{[c]_{x,y}}\Omega_{x,y}$ (resp. at $T_{[c]_x}\Omega_x$) is established
as in the classical case. The
index and nullity of $[c]_{x,y}$ (resp. of $[c]_x$) are defined as usual.
The theory of Bott  concerning
the index and nullity of the $m$-th
iterate of closed geodesics (see  \cite{5, p.495-98} and
\cite{11})
 extends to the case of orbifolds. Anyway we shall show in the next section
that we can consider only the
developable case and, in that case, we can apply the results of
Grove-Tanaka \cite{6}.

\heading 4. Geometric equivalence of closed geodesics \endheading

\subheading{4.1. Geometric equivalence classes}

In this section we describe, for a given closed geodesic of positive length
on an orbifold $Q$, the set of
closed geodesics on
$Q$ which are geometrically equivalent to it, i.e. those closed geodesics
which have the same projection to
$|Q|$. We shall also describe  small invariant tubular neighbourhoods of
their $\Bbb S^1$-orbits. To each
geometric equivalence class of closed geodesics of positive length we
associate a group which is an extension
of $\Bbb Z$ by a finite subgroup. The elements of this group which are not
in this finite subgroup correspond
bijectively  to representative in $\Omega_X$ of elements in this class
which are closed $\Cal G$-geodesic
paths with an initial vector proportional to a given unit vector $\xi$. In
the classical case, this groupis
$\Bbb Z$ and its non zero elements correspond to the multiples of a
primitive closed geodesic.

 Then we prove that the
study of the $\Bbb S^1$-invariant tubular neighbourhoods  in $\Lambda(Q)$
of the closed geodesics in a
geometric equivalence class is equivalent to the corresponding study in the
developable caae.

\mn{\bf 4.1.1. The group attached to a geometric equivalence class.}

Let
$x
\in X$ and let
$( E_x,  e_x)$ be as in 2.4.4 the pointed
principal
$\Cal G$-bundle over $T_xX$ representing the exponential morphism based at
$x$. Let
$\xi \in T_xX$ be a unit vector and let $( E_\xi, e_\xi)$ be
the pointed bundle over $\Bbb R$ pull back of
$(E_x,e_x)$ by the map $t \mapsto t\xi$. To simplify the notations, we note
$p$ the projection $E_\xi \to
\Bbb R$ and $\alpha: E_\xi \to X$ the action map. For a local section $s:
[\tau -\epsilon, \tau +
\epsilon] \to E_\xi$ such that $s(\tau) = e$, we note $\alpha(\dot e)$ the
velocity vector of $\alpha \circ s$
at
$\tau$. The elements
$e \in E_\xi$ with
$p(e)=
\tau
\ne 0$ correspond bijectively  to the equivalence classes of geodesic
$\Cal G$-paths with initial vector equal  to $\tau \xi$. Such an element
represents a closed geodesic
$\Cal G$-path if and only if $\alpha(\dot e) =  \xi$.

The pointed principal $\Cal G$-bundle $(E_\xi,e_\xi)$ is characterized up
to a unique isomorphism
by the following property: For any local section
$s$ of $p$, the composition
$\alpha \circ s$ is a geodesic arc in $X$ with unit speed; moreover
$\alpha(\dot e_\xi) = \xi$.

  Let $e \in E_\xi$ be such that $p(e) = \tau$ and $\alpha(\dot e) = \xi'$.
The pull back of $(E_\xi,e)$ by
the translation
$t
\mapsto t +
\tau$ is canonically isomorphic to $(E_{\xi'},e_{\xi'})$. If $\xi' = \xi$,
this amounts to say that there is a
unique homeomorphism $h: E_\xi \to E_\xi$ projecting to the translation $t
\mapsto t + \tau$, commuting with
the right action of $\Cal G$ and sending $e_\xi$ to $e$. Such a map will be
called an automorphism of $E_\xi$;
it is uniquely determined by $e$. If $c$ is a closed $\Cal G$-geodesic
based at $x$ corresponding to $e
=h(e_\xi)$ and $m$ is a non-zero integer, then $c^m$ corresponds to
$h^m(e_\xi)$.

Let $H$ be the group of such  automorphisms of $E_\xi$. It projects to a
group of translations of
$\Bbb R$ which is discrete because, for $\epsilon$ small enough, there is
no closed $\Cal G$-geodesic loop of
positive length $< \epsilon$ with initial vector proportional to $\xi$.
Assume that this group is not trivial
(equivalently that there is a closed $\Cal G$-geodesic loop of positive
length with initial vector proportional
to
$\xi$); it is an infinite cyclic group generated by an element $\tau_0 >
0$. Let $\phi:H \to \Bbb Z$ be the
homomorphism defined by the relation
$p(h(e_\xi)) = \phi(h)\tau_0$. The kernel of $\phi$ is canonically
isomorphic to $\Cal G_\xi$, the
finite group formed by the elements of $\Cal G$ fixing $\xi$. The action of
$g \in \Cal G_\xi \subset
H$ on an element $e \in E_\xi$ represented by a geodesic $\Cal G$-path $c$
is the element represented by the
path $g.c$.
 We summarize those considerations.

\proclaim{4.1.2. Proposition} Let $\xi\in T_xX$ be a unit vector. Assume
that there is
a closed $\Cal G$-geodesic based at $x$ of positive length with initial vector
proportional to $\xi$. There is an exact sequence of groups
$$1 \to \Cal G_\xi \to H \overset\phi\to\longrightarrow \Bbb Z \to 0$$
where $\Cal G_\xi$ is the group of elements of $\Cal G$ whose differential
fixes
$\xi$. The group $H$ is the group of automorphisms of the principal $\Cal
G$-bundle $E_\xi$ over $\Bbb R$ corresponding to the full geodesic with initial
vector $\xi$.

The elements of $H$ which are not in the subgroup $\Cal G_\xi$ correspond
bijectively to the closed $\Cal G$-geodesics based at $x$ with non-zero
initial vector proportional to
$\xi$.  Let $c$ be a closed $\Cal G$-geodesic path with minimal positive
length $\tau_0$ with initial vector
$\tau_0
\xi$ correponding to an element $h_0 \in H$. Then any closed
$\Cal G$-geodesic of positive length based at $x$ with initial vector
proportional to
$\xi$ corresponding to $h = gh_0^m \in H$, where $g \in
\Cal G_\xi$ and
$m$ is a non-zero integer,  is represented by the $\Cal G$- path $c_h =
g.c^m$.

\endproclaim

\mn{\bf 4.1.3.  Construction of a developable model.} Given a geodesic
$\Cal G$-path
$c=(g_0,c_1,\dots,c_k,g_k)$ from
$x$ to $y$, and a vector $v \in T_xX$, we can construct along each $c_i$ a
parallel vector field $v_i$ such
that $Dg_0(v_1(0)) = v,\  Dg_i(v_{i+1}(t_i)) =v_i(t_i)$. The map $T_xX \to
T_yX$ mapping $v$ to
$Dg_k^{-1}v_k(1)
\in T_yX$ will be called  the parallel transport along $c$.
 It depends only on the equivalence
class of
$c$ and maps the orthogonal to the initial vector of $c$ to the orthogonal
to the terminal vector of $c$. In
particular if $c$ is a closed $\Cal G$-geodesic èpath based at $x$ with
initial vector proportional to $\xi$,
it maps to itself the subspace orthogonal to $\xi$.

Let $N_0$ be the subspace of $T_xX$ orthogonal to $\xi$. Given $h \in H$,
let $\rho(h)$ be the isometry of
$N_0$ given by the inverse of the parallel transport along a closed
geodesic $\Cal G$-loop corresponding to
$h(e_\xi)$. The map $h \mapsto \rho(h)$ is a homomorphism $\rho: H \to
{\text Isom }\  N_0$. Let $N_0^\epsilon$
be the ball of radius $\epsilon$ centered at the origin of $N_0$, and let
$N^\epsilon = \Bbb R \times
N_0^\epsilon$. The group $H$ acts properly on $N^\epsilon$ by the formula
$$h(t,\nu) = (t + \phi(h)\tau_0, \rho(h)\nu).$$
Let $\tilde Q$ be the orbifold $H \backslash N^\epsilon$. We note $\tilde
a$ the map $\Bbb R \to
N^\epsilon=\Bbb R \times N_0^\epsilon$ sending $t$ to $(t,0)$, $\tilde x$
the point $\tilde a(0)$  and $\tilde
\xi$ the velocity vector of
$\tilde a$ at
$t = 0$.

\proclaim{4.1.4. Theorem} For $\epsilon$ small enough, one can define on
$N^\epsilon$ a Riemannian metric
invariant by $H$ and such that $\tilde a$ is a geodesic with speed one, and
an open Riemannian immersion $u$
from the orbifold $\tilde Q = H \backslash N^\epsilon$ to the orbifold $Q =
\Cal G \backslash X$. The
immersion $u$ induces a bijection between the set of geodesic
$H \ltimes N^\epsilon$-loops based at $\tilde x$ with initial vector a
non-zero multiple of $\tilde \xi$ and
 the set of geodesic
$\Cal G$-loops based at $ x$ with initial vector a non-zero multiple of $
\xi$. Moreover the immersion $u$
 induces   isomorphisms of the orbifold neighbourhoods of the corresponding
closed
geodesics and of their $\Bbb S^1$-orbits. \endproclaim

\demo{Proof} Let $h_0$ be an element of $H$ such that $\phi(h_0) = 1$ and
let $c=(g_0,c_1, \dots,
c_k,g_k)$ be a closed geodesic
$\Cal G$-path over a subdivision
$(0=t_0 < t_1 < \dots <t_k=1)$ with initial vector $\tau_0\xi$ whose
equivalence class in $\Omega_x$
corresponds to $h_0$. We can choose $c$ in its equivalence class so that
$g_0 = g_k = 1_x$  and  that there
is an $\epsilon >0$ such that, for each $t \in [t_{i-1},t_i]$, the point
$c_i(t)$ is the center of a convex
geodesic ball $B(c_i(t), 2\epsilon)$ of radius $2\epsilon$. We can also
assume that the length of each $c_i$
is smaller than $\epsilon$. Let $\delta$ be a  positive number smaller than
the numbers $(t_{i-1}
-t_i)/2$. Each $c_i$ can be extended to a geodesic segment $\overline
c_i:]t_{i_1}-\delta,t_i + \delta[ \to X$.

For each $i=1,\dots,k$, let $U_i =]\tau_0(t_{i-1}-\delta),\tau_0(t_i +
\delta)[ \times N^\epsilon_0\subset
 \Bbb R \times N^\epsilon_0 = N^\epsilon$ and let $\tilde X$ be the
disjoint union of the $U_i$, namely the
union of the $\tilde U_i := (\{i\},U_i)$. The composition of the natural
projection $\tilde u_i: \tilde U_i
\to U_i$ with the projection $N^\epsilon \to H \backslash N^\epsilon =
\tilde Q$ gives an atlas of
uniformizing chart for $\tilde Q$. The pseudogroup $\tilde\Cal P$ of change
of charts of this atlas is
generated by the following elements (whenever defined):

i) $\tilde u_i^{-1} \tilde u_{i+1}$ for $i=1,\dots, k-1$,

ii) $\tilde h_0 : = \tilde u_k^{-1} h_0\tilde u_1$,

iii) for each $d \in \Cal G_\xi$ and $i=1,\dots,k$, the maps $\tilde d_i:
\tilde U_i \to \tilde U_i$ sending
$(i,(\tau_0t,\nu))$ to $(i,(\tau_0t,\rho(d)\nu))$.

To describe $\tilde Q$, we  can replace $N^\epsilon$ by $\tilde X$ and the
groupoid $H \ltimes N^\epsilon$ by
the groupoid $\tilde \Cal G$ of germs of $\tilde \Cal P$. The immersion
$\tilde Q \to Q$ will be defined as a
continuous functor $u:\tilde \Cal G \to \Cal G$ that we are going to
construct below.

Given $\nu \in N_0^\epsilon$, let $\overline \nu_i$ be the parallel vector
fields along the $\overline c_i$
such that $\overline \nu_1(0) = \nu$ and $\overline \nu_i(t_i) =
(Dg_i)(\overline \nu_{i+1}(t_i))$. Let
$u_i: \tilde U_i \to X$ be the open embedding mapping $(i, (\tau_0t,\nu))$
to $exp_{ c_i(t)} \overline
\nu_i(t)$. The functor $u$ restricted to the space of units will be the
open immersion  $u: \tilde X \to X$
equal to  $u_i$ on $\tilde U_i$.

We note $\Cal P$ the pseudogroup of local isometries of $X$ whose groupoid
of germs is $\Cal G$. It contains
in particular, for each $i=1,\dots,k-1$, the isometry $\overline g_i:
B(c_{i+1}(t_i),2\epsilon) \to
B(c_i(t_i),2\epsilon)$ whose germ at
$c_{i+1}(t_i)$ is $g_i$  (see 2.1.5). It also contains, for each
$i=1,\dots,k$ and $d \in \Cal G_\xi$, the
isometry
$\overline d_i$ of the ball $B(c_i(t_{i-1}),2\epsilon)$, leaving invariant
$c_i$, the germ at $x$ of
$\overline d_1$ being
$d$
 and the  germ of $\overline d_i$ at $c_i(t_{i-1})$ being the germ of
$\overline g_{i-1}^{-1} \overline d_{i-1}
\overline g_{i-1}$.
We have the following relations, whenever both sides are defined:

i') $u_i(\tilde u_i^{-1} \tilde u_{i+1})u_{i+1}^{-1} = \overline g_i$ for
$i =1,\dots, k-1$,

ii') $u_k\tilde h_0 = u_1$,

iii') $u_i\tilde d_i u_i^{-1} = \overline d_i$ for each $d \in \Cal G_\xi$
and $i=1,\dots,k$.

Therefore we can define the functor $u: \tilde \Cal G \to \Cal G$ by
associating respectively to the germs of
the
$\tilde u_i^{-1}
\tilde u_{i+1}$, of
$\tilde h_0$ and of the $\tilde d_i$, respectively the germs of the
$\overline g_i$, of the identity and of
$\overline d_i$ at the points corresponding to each other by $u: \tilde X
\to X$.

The above relations show that the Riemannian metric on $\tilde X$ which is
the pull back of the Riemannian
metric on $X$ by the immersion $\tilde X \to X$ is $\tilde \Cal
G$-invariant. Hence we obtain on $N^\epsilon$
an $H$-invariant Riemannian metric for which $\tilde a$ is a geodesic.

Instead of looking at the closed geodesic $(H\ltimes N^\epsilon)$-loops
based at $\tilde x$ with initial vector
proportional to $\tilde \xi$, we can equivalently consider the closed
geodesic $\tilde \Cal G$-loops based at
$\tilde x = (1,\tilde x) \in \tilde U_1$ with initial vector proportional
to $\tilde \xi = (1 ,\tilde \xi) \in
T_{\tilde x} \tilde U_1$. For $d \in \Cal G_\xi$, let $\tilde d$ be the
element of $ \tilde \Cal G_{\tilde \xi}$
mapped to
$d$ by $u$.

Let $\tilde c_{h_0} = (\tilde g_0,\tilde c_1,\dots,\tilde c_k,\tilde g_k)$
be the
closed geodesic
$\tilde\Cal G$-path over the subdivision $0=t_0< t_1<\dots <t_k=1$ with
initial vector $\tau_0\tilde\xi$, where
$\tilde c_i: [t_{i-1},t_i] \to \tilde U_i$ is such that $u_i\tilde c_i =
c_i$, $\tilde g_0 = 1_{\tilde
x}$, $\tilde g_k$ is the germ at $\tilde x$ of $\tilde h_0$, and for $i =
1,\dots,k-1$, $\tilde g_i$ is the
germ of $\tilde u_i^{-1}\tilde u_{i+1}$ at $\tilde c_{i+1}$. The functor $u
: \tilde \Cal G \to \Cal G$ maps
$\tilde c_{h_0}$ to $c_{h_0}$, and more generally, for $h=dh_0^m$, where $d
\in \Cal G_\xi$, it maps $\tilde
c_h:=\tilde d.\tilde c_{h_0}^m$ to
$c_h:=d.c_{h_0}^m$. The immersion $u$ induces an immersion $\Omega_u:
\Omega_{\tilde X}(\tilde \Cal G)
\to
\Omega_X(\Cal G)$. It is clear that $\Omega_u$ maps equivariantly  small
neighbourhoods of $[\tilde
c_h]_{\tilde x}$ to small neighbourhoods of $[c_h]_x$ for every $h \in H$.
Moreover the map $\Lambda_u:
\Lambda(\tilde Q) \to \Lambda(Q)$ induced by $\Omega_u$ gives $\Bbb
S^1$-invariant  isomorphisms of  tubular
neighbourhoods of the $\Bbb S^1$-orbits of $[\tilde c_h]$ and $[c_h]$.
$\square$ \enddemo

\mn{\bf 4.1.5.  Description of a $\Bbb S^1$-invariant tubular neighbourhood
of the $\Bbb S^1$-orbit of
$[c_h]$.} To describe such a tubular neighbourhood, using the preceding
theorem, we can as well describe a
$\Bbb S^1$-invariant tubular neighbourhood of the $\Bbb S^1$-orbit of the
corresponding based closed $(H
\ltimes N^\epsilon)$-geodesic $\tilde c_h = (\tilde a^{m\tau_0},h)$, where
$m =\phi(h)$ and $\tilde
a^{m\tau_0}(t) =
\tilde a(m\tau_0t) = (m\tau_0t,0) \in \Bbb R \times N_0^\epsilon  =
N^\epsilon$.

The subgroup of $H$ leaving invariant the $\Bbb R$-orbit of $\tilde c_h$ is
the centralizer $Z_h$ of $h$ in $H$.
The inage of $Z_h$ by $\phi$ is a subgroup of $\Bbb Z$ containing
$\phi(h)=m$, hence it is generated by a
positive integer $d$ dividing  $m$. Comparing with the general situation
described in 3.3.7,
$Y$ corresponds to the connected component of $\Omega_{N^\epsilon}$
containing $[\tilde
c_h]_{\tilde x}$, the group $\overline
\Gamma$ corresponds to $Z_h$ (it is equal to $G$ in our case), and the
homomorphism $\psi: G \to \Bbb Z\  1/r$
to the homomorphism $Z_h \to \Bbb Z\  1/r$, where $r = m/d$, mapping $h'
\in Z_h$ to $\phi(h')/r$. Also the
element
$\gamma_0$ corresponds to
$h \in Z_h$,  and the group $I$ to the quotient $I_h$ of $Z_h$ by the
subgroup generated by $h$.

 Let $\Cal N_0$ be the subspace of the tangent space of
$\Omega_{N^\epsilon}$ at $[\tilde c_h]_{\tilde x}$
orthogonal to the speed vector field along $\tilde c_h$. The elements of
$\Cal N_0$ are the
$H^1$-vector fields
$t
\mapsto v(t)$ along the geodesic
$\tilde{a}^{m\tau_0}$   such that the image of $v(t)$ by the differential
of $h$ is
$v(t+1)$ and
$\int^1_0
\langle v(t),\dot{
\tilde a^{m\tau_0}}(t)
\rangle dt =0$. The group
$Z_h$ acts by isometries on $\Cal N_0$ through  the differential of its
elements acting on $N^\epsilon$,
because those elements commute with $h$. Indeed let us note $v'$ the vector
field  along $\tilde a
^{m\tau_0}$ associating to $t$ the image $v'(t)$ of $v(t-\psi(h'))$ by the
differential $Dh'$ of $h' \in Z_h$.
We have $Dh [v'(t)] = DhDh'[v(t -\psi(h'))] = Dh'[v(t+1-\psi(h'))] =v'(t+1)$.
As  $h$ itself acts trivialy,  this gives an action on $\Cal N_0$ of
the quotient
$I_h$ of
$Z_h$ by the subgroup generated by $h$. This  induces  an action of $I_h$
by isometries on the image $\Cal
D_0$ in
$\Omega_X$ by the exponential map of a small ball in $\Cal N_0$ centered at
$0$.

The homomorphism $\psi: Z_h \to \Bbb Z\  1/r$ induces a homomorphism
$\phi_h:I_h \to \Bbb S^1$ whose image is a
subgroup of order
$r = m/d$. The quotient $\Bbb S^1 \times _{I_h} \Cal D_0$ of $\Bbb S^1
\times \Cal D_0$ by the diagonal action
of
$I_h$ is isomorphic as an orbifold (see 3.3.8) to an $\Bbb S^1$-invariant
tubular neighbourhood
$U$ of $\Bbb S^1.[\tilde c_h]$, or equivalently of  $\Bbb S^1.[c_h]$.

\subheading{4.2. Local homological invariants in a geometric equivalence class}

Following Gromoll-Meyer \cite{5} and Grove-Tanaka \cite{6} we define local
homological invariants for isolated
critical orbits and compare them for geodesics in a geometric equivalence
class.

\mn{4.2.1.\bf Local homological invariants for isolated critical orbits.}
Let $c$ be a closed geodesic $\Cal
G$-path of positive length based at $x \in X$. We have just seen that  a
$\Bbb S^1$-invariant tubular
neighbourhood of the
$\Bbb S^1$-orbit of the equivalence class $[c]$ of $c$ is isomorphic to the
quotient $\Bbb S^1 \times_I \Cal
D_0$ of
$\Bbb S^1
\times \Cal D_0$ by the diagonal action of a finite group $I$ commuting
with the action of $\Bbb S^1$. The
fiber $\Cal D_0$ is isomorphic to the image by the exponential map of a
small ball in the subspace of
$T_{[c]_x}\Omega_X$ orthogonal to the velocity vector field of $c$. We note
$E$ the energy function  on $\Bbb
S^1 \times \Cal D_0$ (it is $\Bbb S^1$-invariant).  We
assume that $\Bbb S^1.[c]$ is an isolated critical orbit and, following
Gromoll and Meyer \cite{5}, we want to
define a local homological invariant $\Cal H_*(\Bbb S^1.[c])$ depending
only of the germ of neighbourhood of
the orbit
$\Bbb S^1.[c]$.

The restriction $E_0$ of the energy function to the fiber $\Cal D_0$ has
$[c]_x$ as an isolated critical point
with critical value $a$. By taking $\Cal D_0$ small enough, we can assume
that $[c]_x$ is the only critical
point of
$E_0$. As in \cite{4}   (see also \cite{5, p.502}), one can construct an
admissible region $W_0$ for the
function
$E_0$; it is constained in $E_0^{-1}[a-\delta,a+\delta]$, where $\delta$ is
a small positive number, and can
be  choosen to be invariant under the action of
$I$.  The local invariant
$\Cal H_*(E_0,[c]_x)$ is defined in \cite{4} as the singular homology
$H_*(W_0,W_0^-)$ with rational
coefficient, where
$W_0^-$ is the intersection of
$W_0$ with $E_0^{-1}(a-\delta)$. Setting $W = \Bbb S^1 \times_IW_0$ and
$W^-=\Bbb S^1 \times_I W_0^-$, the
local homological invariant for the $\Bbb S^1$-orbit of $[c]$ is the finite
dimensionla vector space defined
by
$$\Cal H_*(\Bbb S^1.[c]) : = H_*(W,W^-).$$
It is isomorphic to the subspace $H_*(\Bbb S^1 \times W_0)^I$ of $H_*(\Bbb
S^1 \times W_0)$ of elements left
invariant by the induced action of $I$ on the homology.

The dimension of $\Cal H_k(\Bbb S^1.c)$ is noted $B_k([c])$ and is called
the $k$-th local Betti
number of the isolated orbit $\Bbb S^1.[c]$.

We have the analogue of Lemma 4 in \cite{5} with the same proof. For $b \ge
0$, we note $|\Lambda Q|^b$ the
subspace of $|\Lambda Q|$ which is the inverse image of $[0,b]$ by the
energy function $E: |\Lambda Q| \to
\Bbb R$.

\proclaim{
4.2.2. Lemma} Let $Q$ be a compact Riemannian orbifold. Let $a > 0$ be the
only critical value of the
energy function
$E : |\Lambda Q| \to \Bbb R^+$ in the interval $[a -\epsilon,a+\epsilon]$.
Assume that the critical set  in
$E^{-1}(a)$ consists of finitely many critical orbits $\Bbb S^1.[c^1],\dots
\Bbb S^1.[c^r]$. Then
$$H_*(|\Lambda Q|^{a+\epsilon},|\Lambda Q|^{a-\epsilon}) =  \sum_{i=1}^r
\Cal H_*(\Bbb S^1.[c_i]).$$

In this statement, the space $|\Lambda Q|$ can be replaced by a union of
some of its connected components.
\endproclaim

 We have the analogue of
Corollary 2 of Gromoll-Meyer \cite{5} (see also \cite{6}).

\proclaim{4.2.3. Theorem} Let $\Cal C$ be a geometric equivalence class of
closed geodesics of positive length
on the orbifold $Q =\Cal G \backslash X$. We assume that each $\Bbb
S^1$-orbit in this equivalence class is an
isolated critical manifold for the energy function.

1) There is a constant $B$ such that  $B_k([c]) < B$  for every $k > 2
\text{dim}\ |Q|$ and every
$[c]
\in
\Cal C$.

2) There is a constant $C$ such that  the number of orbits $\Bbb S^1.[c]$,
for $[c] \in \Cal C$, such that
$B_k([c])
\ne 0$ for $k > 2 \text{ dim}\ |Q|$, is  bounded by $C$. \endproclaim

\demo{Proof} Although a more direct proof could be obtained  following the
arguments of Gromoll-Meyer and
Grove-Tanaka, we shall deduce the theorem from the results of Grove-Tanaka
\cite{6}.

Let $\xi  \in T_xX$ be a unit tangent vector such that the elements of
$\Cal C$ are the $\Bbb
S^1$-orbits of elements represented by closed geodesic $\Cal G$-paths with
initial vector proportional to
$\xi$.  Let $N^\epsilon_0$ be an $\epsilon$-ball in the orthogonal to $\xi$
in $T_xX$ centered at $0$. According
to 4.1, there is
 a group $H$ acting properly by isometries on a Riemannian manifold
diffeomorphic to
$ N^\epsilon = \Bbb R
\times  N_0^\epsilon$. The curve
$\tilde a: t \mapsto (t,0)$ is a geodesic line with unit speed left
invariant by $H$.  The elements of $\Cal C$
correspond bijectively to the elements of the geometric class $\tilde \Cal
C$ of closed geodesics on the
orbifold
$\tilde Q = H\backslash N^\epsilon$ having the same image in $|\tilde Q|$
as $\tilde a$ and the corresponding
elements have the same local homological invariants. Therefore it is
sufficient to prove the theorem for
$\tilde \Cal C$.

We have a surjective homomorphism $\phi: H \to \Bbb Z$ ; the kernel of
$\phi$ is the finite subgroup  of $H$ fixing the velocity vector
$\dot{\tilde a}(0) =\tilde \xi$. Recall that
the
$\Bbb S^1$-orbits of the elements of
$\Cal C$ are parametrized by the conjugacy classes of the elements of $H$
whose image by $\phi$ is non zero.
To each element $h$ of $H$ is associated the element of $\tilde \Cal C$
represented by the geodesic path
$\tilde c_h =(\tilde a^{\phi(h)\tau_0},h)$, where $\tau_0$ is the minimal
length of the closed geodesics in the
class
$\Cal C$.

As $H$ is an extension of $\Bbb Z$ by a finite group, we can find an
element $h_1 \in H$ in the center of
$H$ with $\phi(h_1) = n_1 >0$. Let $\Gamma$ be  the quotient of $H$ by the
subgroup
$H_1$ generated by
$h_1$. It is a finite group acting on the quotient Riemannian manifold $M:=
H_1 \backslash
N^\epsilon$. Let $\pi:H \to \Gamma$ be the quotient projection. The
natural projection
$ N^\epsilon
\to M$ induces an isomorphism
$\tilde Q=H\backslash
 N^\epsilon \to \overline Q:=\Gamma\backslash M$ of the quotient orbifolds.
The composition of
$\tilde a:\Bbb R
\to
 N^\epsilon$ with the natural projection  onto $M$ is a closed geodesic
$\overline a:\Bbb R \to M$.
Therefore we can replace $\tilde \Cal C$ by the geometric equivalence class
$\overline {\Cal C}$ of closed
geodesics on $\overline Q$ which have the same image in $|\overline Q|$ as
the image of $\overline a$. To each
$h \in H$ is associated the closed based geodesic $(\Gamma \ltimes M)$-path
with initial vector
proportional to $\overline \xi:= \dot{\overline a}(0) $  represented by
$\overline c_h:= (\overline a^{\phi(h)\tau_0}, \pi(h))$.

For each $\gamma \in \Gamma $, let $\Cal C_{\gamma}$ be the set of
$\gamma$-invariant geodesics (in the sense of Grove) in $M$ with non-zero
initial vector proportional to
$\dot{\overline a}(0)$. These geodesics correspond bijectively to the
closed geodesics $\overline c_h$ with
$\pi(h) = \gamma$ and $\phi(h) \ne 0$. As $\overline {\Cal C}$ is the
finite union of the $\Bbb S^1$-orbits
of the elements of $\Cal C_\gamma$, it is sufficient to prove the analogue
of the theorem where $\Cal
C$ is replaced by $\Cal C_\gamma$.

In the notations of Grove-Tanaka, let $\Lambda(M,\gamma)$ be the space of
$\gamma$-invariant
$H^1$-curves on $M$. (In our notations this corresponds to a union of
connected components of
$\Omega_M(\Gamma \ltimes M)$.). The natural  action of $\Bbb R$ on
$\Lambda(M,\gamma)$ gives the
action of the circle
$\Bbb R/s\Bbb Z$ considered in \cite{6}, where $s$ is the order of
$\gamma$. The restriction of the energy
function to
$\Lambda(M,\gamma)$ is noted $E^{\gamma}$. Let $\mu: \Lambda(M,\gamma) \to
|\Lambda \overline Q|$ be the
composition of the inclusion  $\Lambda(M,\gamma) \to \Omega_M(\Gamma
\ltimes M)$ with the quotient map to
$\Gamma\backslash \Omega_M(\Gamma \ltimes M) =|\Lambda \overline Q|$. The
image of $\mu$ is an open set
isomorphic to the quotient of $\Lambda(M,\gamma)$ by the action of the
centralizer $Z_\gamma$ of $\gamma$ in
$\Gamma$. The map $\mu$ sends  $\Bbb R/s \Bbb Z$-orbits to  $\Bbb
S^1$-orbits. In a small $\Bbb
S^1$-invariant tubular neighbourhood of $[\overline c_h]$, where $h \in
\pi^{-1}(\gamma)$, following
Gromoll-Meyer \cite{5},  we have constructed in 4.2.1 an admissible region
$W$ and defined $B_k([\overline
c_h])$ as the dimension of $H_k(W,W^-)$. The inverse image $\mu^{-1}(W)$ of
$W$ is a $\Bbb R/s \Bbb
Z$-invariant admissible region  and $ \text{dim }
H_k(\mu^{-1}(W),\mu^{-1}(W^-))$ is the
local Betti number defined in Grove-Tanaka \cite{6, p. 44} noted there
$B_k(\overline
a^{\phi(h)\tau_0},\gamma)$. As
$W = Z_\gamma
\backslash \mu{^{-1}}(W)$, the vector space $H_k(W,W^-)$ is isomorphic to
the vector subspace of
$H_k(\mu^{-1}(W),\mu^{-1}(W^-))$ left invariant by $Z_\gamma$. Therefore
$B_k([\overline c_h]) \le
B_k(\overline
a^{\phi(h)\tau_0},\gamma)$, and we can use the results of \cite{6} where
the existence of constants like $B$
and $C$ are established for $B_k(\overline c_h,\gamma)$, where $\overline
c_h \in \Cal C_\gamma$. $\square$
\enddemo

\heading 5. Existence of geodesics \endheading

\subheading{5.1. Existence of at least one  closed geodesic of positive
length}

\proclaim{5.1.1. Theorem} Let $Q$ be a compact connected Riemannian
orbifold. There exists at least
one closed geodesic on $Q$ of positive length in the following cases:

a) $Q$ is not developable,

b) the fundamental group of $Q$ has an element of infinite order or is finite.
\endproclaim

\demo{Proof} a)  If $Q = \Cal G\backslash X$ is not developable, there is a
point $x \in X$ and
a non trivial element $g$ in the group $ \Cal G_x$ of elements of $\Cal G$
fixing $x$ such that
the closed loop based at $x$ represented by the  $c=(1_x,c_1,g)$, where
$c_1:[0,1] \to X$ is the
constant map to $x$, is homotopically trivial (see for instance \cite{2}).
Hence in
$|\Lambda Q|$ there is a continuous path joining the point $z$ of
$|\Lambda^0Q|$ represented by
$c$ to a point
$z'$ of
$|\Lambda^0Q|$ represented by a constant loop. Those two points are in
distinct components of
$|{\Lambda^0}Q|$, and the components are all compact. Moreover the sets of
point
$|{\Lambda}^\epsilon Q|$ of $|\Lambda Q|$ for which the energy function is
smaller than $\epsilon$
for various $\epsilon >0$ form a fundamental system of neighbourhoods of
$|{\Lambda}^0Q|$ (to see this adapt the proof given in \cite{11, pages
30-31}, using 2.4.6).
Therefore the family of paths in $|\Lambda Q|$ joining $z$ to $z'$ is a
$\phi$-family $\mod
|{\Lambda}^0Q|$ and we can apply 3.4.6.

b) The energy function restricted to the  connected component of
$|\Lambda Q|$ corresponding to
an element of infinite order of the fundamental group of $Q$ attains its
infimum (cf. 3.4.7) at
some point; this point is necessarily of  positive length, as this curve
represents an element of
infinite order. In the case where the fundamental group of
$Q$ is finite, then either $Q$ is not developable and we can apply a), or
its universal covering is
a compact Riemannian manifold $M$. Then   the classical result of Fet
implies the existence of a
closed geodesic on
$M$ of positive length and its projection  gives a closed geodesic on $Q$
of positive length.
$\square$ \enddemo

\remark{5.1.2. Remark} For the existence of a closed geodesic of positive
length on
compact orbifolds, the only case left open by the preceding theorem would
be the following one. Let
$\Gamma$ be an infinite group all of whose elements are of finite order
acting properly on a
Riemannian manifold $M$ by isometries with compact quotient $Q = \Gamma
\backslash M$; does
there exists a non constant geodesic $c:[0,1] \to M$ and an element $\gamma
\in \Gamma$ such
that the differential of $\gamma$ maps $\dot c(0)$ to $\dot c(1)$ ?. Note
that such a group
$\Gamma$ would be finitely presented, and no examples of such groups are
known yet. \endremark

\remark{5.1.3. Remark} The argument used in a) applies whenever  two
distinct connected
components of $|{\Lambda}^0Q|$ are contained in the same  connected
component of $|\Lambda
Q|$.

For instance, consider on $\Bbb S^2$ a Riemannian metric invariant by a
rotation $\rho$ of order $n$ fixing
the north pole $N$ and the south pole $S$. The quotient of $\Bbb S^2$ by
the group generated by $\rho$ is a
Riemannian orbifold. Let $k$ be an integer not divisible by $n$; the
elements of $|\Lambda Q|$ represented
by the closed
$\Cal G$-paths $(c,\rho^k)$, where $c:[0,1] \to \Bbb S^2$ is either the
constant map to $N$ or to $S$ are
in the same connected component of $|\Lambda Q|$, but in distinct
components of $|{\Lambda}^0Q|$.
Therefore the argument in a) shows that there exists a geodesic segment
$c:[0,1]
\to
\Bbb S^2$ of positive length such that the differential of $\rho^k$ maps
$\dot c(0)$ to $\dot c(1)$.
\endremark

\subheading{4.2. Existence of infinitely many geodesics}

The following is the generalization of a theorem of Serre \cite{ 16}.

\proclaim{4.2.1. Theorem} Let $Q$ be a compact connected orbifold. Given
two points
$\overline x$ and
$\overline x'$ of $|Q|$, there exist an infinity of geodesics from
$\overline x$ to $\overline x'$.
\endproclaim

\demo{Proof} Let $Q =\Cal G \backslash X$ and let $x$ and $x'$ be two
points of $X$ projecting
to
$\overline x$ and $\overline x'$. We have to prove that there exists an
infinity of equivalence
classes of geodesic $\Cal G$-paths from $x$ to $x'$. The connected
components of
$\Omega_{x,x'}$ correspond bijectively to the elements of the fundamental
group of $Q$ (i.e. of
$BQ$, see 3.3.5). The energy function assumes its infimum on each connected
component (see 3.4.7).
Therefore it suffices to consider the case where the fundamental group of
$Q$ is finite and, after
passing to the universal covering,  the case where $Q$ is simply connected.
Then the fundamental
rational class of $Q$ (see 2.5.3) gives a non-trivial element in
$H_n(BQ,\Bbb Q)$, and
by Serre \cite{16}, p. 484,
the  rational Betti numbers $b_i$ of $\Omega_{x,x'}$, which has the same
homotopy type as the loop
space of $BQ$, do not vanish for an infinity of value of $i$. This would
contradict the existence of only
finitely many critical points of the energy function on $\Omega_{x,y}$
(see Gromoll-Meyer \cite{4} or the
more direct arguments of Seifert-Threlfall \cite {15}).
$\square$ \enddemo

 We next prove the analogue of Gromoll-Meyer theorem (see also
theorem 4.1 of Grove-Tanaka \cite{6}) about the existence of infinitely
many geometrically distinct closed
geodesics.

\proclaim{5.2.2. Theorem} Let $Q= \Cal G\backslash X$ be a compact connected
Riemannian orbifold. Assume there is only a finite number of geometric
equivalence
classes of closed geodesics on $Q$ of positive length. There is a constant $D$
such that the rational Betti numbers $B_k$ of $|\Lambda Q|$ are bounded by
$D$ for
$k \ge 2\text{ dim } Q$. \endproclaim

\demo{Proof} As the proof of Theorem 4 of Gromoll-Meyer \cite{5}, or
Theorem 4.1 in
Grove-Tanaka \cite{6} using 4.2.3. $\square$ \enddemo

The next corollary follows from 5.2.2 and the Vigu\'e-Sullivan theorem
\cite{18}.

\proclaim{5.2.3. Corollary} Let $Q$ be a compact simply connected Riemannian
orbifold. Assume that the rational cohomology of $|Q|$ is not generated by a
single element. Then there are  infinitely many geometric equivalence
classes of
closed geodesics of positive length on $Q$. \endproclaim

\demo{Proof} As $BQ$ is simply connected, $|Q|$ is also simply
connected (the natural map $\pi_1(BQ) \to \pi_1(|Q|)$ is always surjective).
Moreover the projection $BQ \to |Q|$ induces an isomorphism on the rational
homology, so it is a rational homotopy equivalence. This implies that the map
$\Lambda BQ \to \Lambda|Q|$ induced on the free loop spaces is also a rational
homotopy equivalence. So $\Lambda|Q|$ has the same rational Betti numbers
as $\Lambda BQ$

On the other hand  the projection $B\Lambda Q \to |\Lambda Q|$ induces an
isomorphism in rational cohomology by 2.5.3. As $B\Lambda^c Q$ is weakly
homotopy equivalent to
$\Lambda BQ$ by 3.2.2 and to
$B(\Lambda Q)$ by 3.3.5, it follows  finally that $\Lambda |Q|$ has the
same rational Betti numbers as $|\Lambda Q|$.

As $|Q|$ is simply connected, the rational Betti numbers of its free loop space
$\Lambda |Q|$ are all finite by Serre \cite{16}. If there would be only a
finite
number of geometric equivalence classes of closed geodesics on $Q$ of positive
length, then the rational Betti numbers of $|\Lambda Q|$ would be uniformly
bounded by 5.2.2, so also those of $\Lambda |Q|)$. By \cite{18} this would
imply that
the rational cohomology of $|Q|$ would be generated by a single element,
contrary
to the hypothesis. $\square$ \enddemo

The next theorem is the extension to the case of orbifolds of a result of
Viktor Bangert and Nancy Hingston
\cite{1}.

\proclaim{5.2.4. Theorem} Let $Q$ be a connected compact Riemannian
orbifold of dimension $> 1$ whose
fundamental group is infinite abelian. Then there is an infinity of
geometrically distinct closed
geodesics
 of positive length on $Q$. \endproclaim

\demo{Proof} As mentioned in \cite{1}, the crucial case is when the
fundamental group of $Q$
is infinite cyclic. We can also assume $Q$ orientable. We follow closely the
arguments of
\cite{1}; as far as the homotopy properties are concerned,  the role of
 the compact Riemannian manifold $M$ is played by $BQ$. The fundamental
class of $Q$ gives a
nontrivial element of the homology with rational coefficients of $BQ$ in
dimension $>1$ (see 2.5.3), hence
$BQ$ has not the homotopy type of a circle and $\pi_n(BQ) \neq 0$ for some
minimal $n >1$. Let
$z \in BQ$ be a base point. The connected components of $\Omega_zBQ$
correspond to the elements
of the fundamental group of $Q$. For each integer
$r>0$ choose a loop
$l_r$ at
$z$ whose homotopy class if the  $r$-th power $t^r$ of
 a generator $t$ of the fundamental group of
$Q$.

The lemmas 1 and 2 of \cite{1} provide the existence  of a positive integer
$k$ such that, for
all positive integer $m$, there exist elements $\alpha_m \in
\pi_{n-1}(\Omega_zBQ, l_{mk})$ such
that the image of $\alpha_m$ in $\pi_{n-1}(\Lambda BQ, l_{mk})$ is non
trivial. The element
$\alpha_m$ is described as follows. There is a map $f: (\Bbb S^{n-1},*) \to
(\Omega_zBQ, z) $
(where $z$ is identified to the constant loop at $z$) so that  the element
$\alpha_m$ is represented by the
map $f_m$ obtained by composing with $l_{mk}$ the loops in the image of $f$.

 By 3.2.2 and 3.3.5, $\Omega_zBQ$ has the same weak homotopy type as
$\Omega_x$, where $x\in X$ and $z$
have the same projection to $|Q|$. Let
$[c_m]_x
\in \Omega_x$ be a $\Cal G$-loop based at $x$ of minimal energy in the
homotopy class of $t^{mk}$. Let
$f': (\Bbb S^{n-1},*) \to (\Omega_x, [c_0]_x)$ be a map whose homotopy
class corresponds to the
homotopy class of $f$. Let $f'_{m}: (\Bbb S^{n-1},*) \to
(\Omega_x,[c_{m}]_x)$ mapping $y \in
\Bbb S^{n-1}$ to the composition of $f'(y)$ with $[c_{m}]_x$. Considered as
a morphism from $\Bbb
S^{n-1}$ to the orbifold $\Lambda Q$, it is not homotopic to a constant.

Let $\kappa_m$ be the infimum of the energy of the $\Cal G$-loops in the
homotopy class $t^{km}$.
Let
$F_m \subset |\Lambda Q|$ be the image of the composition of  $f'_m$  with
the projection
to
$|\Lambda Q|$ and let $\tau_m$ be the critical value of the $\phi$-family
$\Cal F$ (see 3.4.5)
whose elements are the subsets
$\phi_t(F_m)$ ($t>0$), where $\phi_t$ is the flow of $-\text{grad  }E$.
Following \cite{1} one can
assume that $\kappa_m > \tau_m$, because otherwise either there would exist
a continuous family of
geometrically distinct  closed geodesics with energy $\kappa_m$, or there
would exist in the homotopy
class of $f'_m$, by 3.4.4, a morphism whose image in $|\Lambda Q|$ would be
contained in a $\Bbb
S^1$-orbit; but this would contradict the condition that $f'_m$ is not
homotopic in $\Lambda Q$ to a
constant. The final part of the proof is as in
\cite{1}.
$\square$

\medskip

To end up we reformulate in our framework some  results of Grove and
Halperin \cite{7}.

\proclaim{5.2.5. Theorem of Grove and Halperin} Let $Q$ be an orbifold
quotient of a closed
simply connected Riemannian manifold
$X$ by a finite group of isometries.

1) If the dimension of $Q$ is odd, there is a closed geodesic of positive
length on $Q$ in each
connected component of
$\Lambda Q$.

2) Assume that the rational homotopy groups $\pi_k(X)$ of $X$ are non zero
for an infinity
of $k$. Then in each component of $\Lambda Q$, there exist an infinity of
geometrically
distinct closed geodesics on $Q$ of positive length. \endproclaim

\demo{Proof} This theorem follows from Therorem A and B in \cite{7, p. 173}
using the remarks
at the end of 3.3.6. $\square$ \enddemo

\enddemo

\enddocument

\bye